\title{Additive Approximation for Edge-Deletion Problems
}
\author{Noga Alon \thanks{
Schools of Mathematics and Computer Science, Raymond and Beverly
Sackler Faculty of Exact Sciences, Tel Aviv University, Tel Aviv
69978, Israel and IAS, Princeton, NJ 08540, USA. Email:
nogaa@tau.ac.il. Research supported in part by the Israel Science
Foundation, by a USA-Israeli BSF grant,
by the Hermann Minkowski Minerva Center for Geometry
at Tel Aviv University and by the Von Neumann Fund. } \and Asaf
Shapira
\thanks{School of Computer Science,
Raymond and Beverly Sackler Faculty of Exact Sciences, Tel Aviv
University, Tel Aviv, Israel. Email: asafico@tau.ac.il. This work
forms part of the author's Ph.D. thesis. Research supported in
part by a Charles Clore Foundation Fellowship and an IBM Ph.D.
Fellowship.} \and Benny Sudakov
\thanks{ Department of Mathematics, Princeton University,
Princeton, NJ 08544, USA. E-mail: bsudakov@math.princeton.edu.
Research supported in part by NSF CAREER award DMS-0546523, NSF grant
DMS-0355497, USA-Israeli BSF grant and by an Alfred P. Sloan fellowship.
} }

\date{}

\documentclass [letterpaper,11pt]{article}
\usepackage{amsfonts}
\usepackage{amsmath}

\newcommand{\qed}{\hspace*{\fill} \rule{7pt}{7pt}}

\newcommand{\ignore}[1]{}

\newtheorem{theo}{Theorem}[section]
\newtheorem{lemma}{Lemma}[section]
\newtheorem{coro}[lemma]{Corollary}

\newtheorem{definition}[lemma]{Definition}
\newtheorem{prop}[lemma]{Proposition}
\newtheorem{comment}[lemma]{Comment}
\newtheorem{claim}[lemma]{Claim}

\begin{document}
\maketitle
\begin{abstract}

A graph property is {\em monotone} if it is closed under removal
of vertices and edges. In this paper we consider the following
algorithmic problem, called the
edge-deletion problem; given a monotone property ${\cal P}$ and a
graph $G$, compute the smallest number of edge deletions that are
needed in order to turn $G$ into a graph satisfying ${\cal P}$. We
denote this quantity by $E'_{{\cal P}}(G)$. The first result of
this paper states that the edge-deletion problem can be
efficiently approximated for any monotone property.

\begin{itemize}
\item For any fixed $\epsilon > 0$ and any monotone property ${\cal
P}$, there is a deterministic algorithm, which given a graph $G=(V,E)$
of size $n$, approximates $E'_{{\cal P}}(G)$ in linear time $O(|V|+|E|)$ to
within an additive error of $\epsilon n^2$.
\end{itemize}

Given the above, a natural question is for which monotone
properties one can obtain better additive approximations of
$E'_{{\cal P}}$. Our second main result essentially resolves this
problem by giving a precise characterization of the monotone graph
properties for which such approximations exist.

\begin{enumerate}
\item[(1)] If there is a bipartite graph that does not satisfy
${\cal P}$, then there is a $\delta > 0$ for which it is possible
to approximate $E'_{{\cal P}}$ to within an additive error of
$n^{2-\delta}$ in polynomial time.

\item[(2)] On the other hand, if all bipartite graphs satisfy ${\cal P}$,
then for any $\delta > 0$ it is $NP$-hard to approximate $E'_{{\cal
P}}$ to within an additive error of $n^{2-\delta}$.
\end{enumerate}
While the proof of (1) is relatively simple, the proof of (2)
requires several new ideas and involves tools from Extremal Graph
Theory together with spectral techniques.
Interestingly,
prior to this work it was not even known that
computing $E'_{{\cal P}}$ {\em precisely}
for the properties in (2)
is $NP$-hard. We thus answer (in a strong form) a question of
Yannakakis, who asked in 1981 if it is possible to find a
large and natural family of graph properties for which computing
$E'_{{\cal P}}$ is $NP$-hard.

\end{abstract}

\newpage

\section{Introduction}\label{intro}

\subsection{Definitions, background and motivation}

The topic of this paper is graph modification problems,
namely problems of the type: "given a graph $G$, find the smallest
number of modifications that are needed in order to turn $G$ into
a graph satisfying property ${\cal P}$". The main two types of
such problems are the following, in {\em node modification}
problems, one tries to find the smallest set of vertices, whose
removal turns $G$ into a graph satisfying ${\cal P}$, while in
{\em edge modification} problems, one tries to find the smallest
number of edge deletions/additions that turn $G$ into a graph
satisfying ${\cal P}$. In this paper we will focus on
edge-modification problems. Before continuing with the
introduction we need to introduce some notations.

For a graph property ${\cal P}$, let ${\cal P}_n$ denote the set
of graphs on $n$ vertices which satisfy ${\cal P}$. Given two
graphs on $n$ vertices, $G$ and $G'$, we denote by $\Delta(G,G')$
the edit distance between $G$ and $G'$, namely the smallest number
of edge additions and/or deletions that are needed in order to
turn $G$ into $G'$. For a given property ${\cal P}$, we want to
denote how far is a graph $G$ from satisfying ${\cal P}$. For
notational reasons it will be more convenient to normalize this
measure so that it is always in the interval $[0,1]$ (actually
$[0,\frac12]$). We thus define

\begin{definition}\noindent{\bf ($E_{{\cal P}}(G)$)}\label{DefEdit}
For a graph property ${\cal P}$ and a graph $G$ on $n$ vertices,
let
$$
E_{{\cal P}}(G) = \min_{G' \in {\cal P}_n}
\frac{\Delta(G,G')}{n^2}\;.
$$
\end{definition}
In words, $E_{{\cal P}}(G)$ is the minimum edit distance of $G$ to
a graph satisfying ${\cal P}$ after normalizing it by a factor of
$n^2$.

Graph modification problems are well studied
computational problems. In 1979, Garey and
Johnson \cite{GJ} mentioned 18 types of vertex and edge
modification problems. Graph modification problems were
extensively studied as these problems have applications in several
fields, including Molecular Biology and Numerical Algebra. In
these applications a graph is used to model experimental data,
where edge modifications correspond to correcting errors in the
data: Adding an edge means correcting a false negative, while
deleting an edge means correcting a false positive. Computing
$E_{{\cal P}}(G)$ for appropriately defined properties ${\cal P}$
have important applications in physical mapping of DNA (see
\cite{CMNR}, \cite{GGKS} and \cite{GKS}). Computing $E_{{\cal
P}}(G)$ for other properties arises when optimizing the running
time of performing Gaussian elimination on a sparse symmetric
positive-definite matrix (see \cite{Rose}). Other modification
problems arise as subroutines for heuristic algorithms for
computing the largest clique in a graph (see \cite{X}). Some edge
modification problems also arise naturally in optimization of
circuit design \cite{EC}. We briefly mention that there are also
many results about {\em vertex} modification problems, notably
that of Lewis and Yannakakis \cite{LY}, who proved that for any
nontrivial hereditary property ${\cal P}$, it is $NP$-hard to compute
the smallest number of vertex deletions that turn a graph into
one satisfying ${\cal P}$. (A graph property is hereditary if
it is closed under removal of vertices.)

A graph property is said to be monotone if it is closed under
removal of both vertices and edges. Examples of well studied
monotone properties are $k$-colorability, and the property of being
$H$-free for some fixed graph $H$. (A graph is $H$-free if it contains
no copy of $H$ as a not necessarily induced subgraph.)
Note, that when trying to turn a graph into one satisfying a
monotone property we will only use edge deletions. Therefore, in
these cases the problem is sometimes called {\em edge-deletion}
problem. Our main results, presented in the following subsections,
give a nearly complete answer to the hardness of additive
approximations of the edge-deletion problem for monotone properties.

\subsection{An algorithm for any monotone property}

Our first main result in this paper states that for any graph
property ${\cal P}$ that belongs to the large, natural and well
studied family of monotone graph properties, it is possible to
derive efficient approximations of $E_{{\cal P}}$.

\begin{theo}\label{t1}
For any fixed $\epsilon > 0$ and any monotone property ${\cal P}$ there
is a {\bf deterministic} algorithm that given a graph $G$ on $n$
vertices computes
in time $O(n^2)$ a real $E$ satisfying $|E-E_{{\cal P}}(G)| \leq
\epsilon$.
\end{theo}
Note, that the running time of our algorithm is of type
$f(\epsilon)n^2$, and can in fact be improved to linear in the size of the input
by first counting the number of edges,
taking $E=0$ in case the graph has less than $\epsilon n^2$ edges.
We note that Theorem \ref{t1} was not known for
many monotone properties. In particular, such an approximation
algorithm was not even known for the property of being triangle-free
and more generally for the property of being $H$-free for any
non-bipartite $H$.

Theorem \ref{t1} is obtained via a novel structural graph theoretic
technique. One of the
applications of this technique (roughly) yields that every graph
$G$, can be approximated by a small weighted graph $W$, in such a
way that $E_{{\cal P}}(G)$ is approximately the optimal solution of
a certain related problem (explained precisely in Section
\ref{overview}) that we solve on $W$. The main usage of this new
structural-technique in this paper is in proving Lemmas \ref{Easy}
and \ref{Big} that lie at the core of the proof of Theorem
\ref{t1}. This new technique, which may very well have other
algorithmic and graph-theoretic applications, applies a result of
Alon, Fischer, Krivelevich and Szegedy \cite{AFKS} which is a
strengthening of Szemer\'edi's Regularity Lemma \cite{Sz}. We then
use an efficient algorithmic version of the regularity lemma, which
also implies an efficient algorithmic version of the result of
\cite{AFKS}, in order to transform the existential structural result
into the algorithm stated in Theorem \ref{t1}.

We further use our structural result in order to prove the
following concentration-type result regarding the edit distance of
subgraphs of a graph.

\begin{theo}\label{t2}
For every $\epsilon$ and any monotone property ${\cal P}$ there is a
$d=d(\epsilon,{\cal P})$ with the following property: Let $G$ be any
graph and suppose we randomly pick a subset $D$, of $d$ vertices
from $V(G)$. Denote by $G'$ the graph induced by $G$ on $D$. Then,
$$
Prob[~|E_{{\cal P}}(G')-E_{{\cal P}}(G)| > \epsilon] < \epsilon\;.
$$
\end{theo}

\bigskip

An immediate implication of the above theorem is the following,

\begin{coro}\label{t3}
For every $\epsilon > 0$ and any monotone property ${\cal P}$ there
is a {\bf randomized} algorithm that given a graph $G$ computes in
time $O(1)$ a real $E$ satisfying $|E-E_{{\cal P}}(G)| \leq
\epsilon$ with probability at least $1-\epsilon$.
\end{coro}

We stress that there are some computational subtleties regrading the
implementation of the algorithmic results discussed above. Roughly
speaking, one should define how the property ${\cal P}$ is "given"
to the algorithm and also whether $\epsilon$ is a fixed constant or
part of the input. These issues are discussed in Section
\ref{SecAlgo}.

It is natural to ask if the above results can be extended to the larger
family of hereditary properties, namely, properties closed under removal
of vertices, but not necessarily
under removal of edges. Many natural properties such as being
Perfect, Chordal and Interval are hereditary non-monotone
properties. By combining the ideas we used in order to prove
Theorem \ref{t1} along with the main ideas of \cite{ASher} it can
be shown that Theorem \ref{t1} (as well as Theorem \ref{t2} and
Corollary \ref{t3}) also hold for any hereditary graph property.

\subsection{On the possibility of better approximations}

Theorem \ref{t1} implies that it is possible to efficiently
approximate the distance of an $n$ vertex graph from any monotone
graph property ${\cal P}$, to within an error of $\epsilon n^2$ for
any $\epsilon > 0$. A natural question one can ask is for which
monotone properties it is possible to improve the additive error to
$n^{2-\delta}$ for some fixed $\delta > 0$. In the terminology of
Definition \ref{DefEdit}, this means to approximate $E_{{\cal P}}$
to within an additive error of $n^{-\delta}$ for some $\delta
> 0$. Our second main result in this paper is a precise
characterization of the monotone graph properties for which such a
$\delta > 0$ exists\footnote{We assume henceforth that ${\cal P}$
is not satisfied by all graphs.}.

\begin{theo}\label{t4}
Let ${\cal P}$ be a monotone graph property. Then,
\begin{enumerate}
\item If there is a bipartite graph that does not satisfy
${\cal P}$, then there is a fixed $\delta > 0$ for which it is
possible to approximate $E_{{\cal P}}$ to within an additive error
of $n^{-\delta}$ in polynomial time.

\item On the other hand, if all bipartite graphs satisfy ${\cal P}$,
then for any fixed $\delta > 0$ it is $NP$-hard to approximate $E_{{\cal
P}}$ to within an additive error of $n^{-\delta}$.
\end{enumerate}
\end{theo}

While the first part of the above theorem follows easily from
the known results about the Tur\'an numbers of bipartite graphs
(see, e.g., \cite{W}), the proof of the second item involves various
combinatorial tools. These include Szemer\'edi's Regularity Lemma,
and a new result in Extremal Graph Theory, which is stated in
Theorem \ref{t81} (see Section \ref{thm2}) that extends the main
result of \cite{BSTT} and \cite{BKS}. We also use the basic approach
of \cite{Al1}, which applies spectral techniques to obtain an
$NP$-hardness result by embedding a blow-up of a sparse instance to
a problem, in an appropriate dense pseudo-random graph. Theorem
\ref{t81} and the proof technique of Theorem \ref{t4} may be useful
for other applications in graph theory and in proving hardness
results. As in the case of Theorem \ref{t1}, the second part of
Theorem \ref{t4} was not known for many specific monotone properties. For
example, prior to this paper it was
not even known that it is $NP$-hard to {\em precisely} compute
$E_{{\cal P}}$, where ${\cal P}$ is the property of being
triangle-free. More generally, such a result was not known for the
property of being $H$-free for any non-bipartite $H$.

\subsection{Related work}
Our main results form a natural continuation and extension of
several research paths that have been extensively studied. Below
we survey some of them.

\subsubsection{Approximations of graph-modification problems}

As we have previously mentioned many practical optimization problems
in various research areas can be posed as the problem of computing
the edit-distance of a certain graph from satisfying a certain
property. Cai \cite{C} has shown that for any hereditary property,
which is expressible by a finite number of forbidden induced
subgraphs, the problem of computing the edit distance is
fixed-parameter tractable. Khot and Raman \cite{KR} proved that for
some hereditary properties ${\cal P}$, finding in a given graph $G$,
a subgraph that satisfies ${\cal P}$ is fixed-parameter tractable,
while for other properties finding such a subgraph is hard in an
appropriate sense (see \cite{KR}).

Note that Theorem \ref{t1} implies that if the edit distance (in
our case, number of edge removals) of a graph from a property
is $\Omega(n^2)$, then it can be approximated to within any
{\em multiplicative} constant $1+\epsilon$.

\subsubsection{Hardness of edge-modification problems} Natanzon,
Shamir and Sharan \cite{NSS} proved that for various hereditary
properties, such as being Perfect and Comparability, computing
$E_{{\cal P}}$ is $NP$-hard and sometimes even $NP$-hard to
approximate to within some constant. Yannakakis \cite{Y} has shown
that for several graph properties such as outerplanar, transitively
orientable, and line-invertible, computing $E_{{\cal P}}$ is
$NP$-hard. Asano \cite{Asano} and Asano and Hirata \cite{AH} have
shown that properties expressible in terms of certain families of
forbidden minors or topological minors are $NP$-hard.

The $NP$-completeness proofs obtained by Yannakakis in \cite{Y},
were add-hoc arguments that applied only to specific properties.
Yannakakis posed in \cite{Y} as an open problem, the possibility of
proving a general $NP$-hardness result for computing $E_{{\cal P}}$
that will apply to a general family of graph properties. Theorem
\ref{t4} achieves such a result even for the seemingly easier
problem of approximating $E_{{\cal P}}$.

\subsubsection{Approximation schemes for "dense" instances}
Fernandez de la Vega \cite{FDLV} and Arora, Karger and Karpinski
\cite{AKK} showed that many of the classical $NP$-complete problems
such as MAX-CUT and MAX-3-CNF have a PTAS when the instance is
dense, namely if the graph has $\Omega(n^2)$ edges or the 3-CNF
formula has $\Omega(n^3)$ clauses. Approximations for dense
instances of Quadratic Assignment Problems, as well as for
additional problems, were obtained by Arora, Frieze and Kaplan
\cite{AFK}. Frieze and Kannan \cite{FK2} obtained approximations
schemes for several dense graph theoretic problems via certain
matrix approximations. Alon, Fernandez de la Vega, Kannan and
Karpinski \cite{ADKK} obtained results analogous to ours for any
dense Constraint-Satisfaction-Problem via certain sampling
techniques. It should be noted that all the above approximation
schemes are obtained in a way similar to ours, that is, by first
proving an {\em additive} approximation, and then arguing that in
case the optimal solution is large (that is, $\Omega(n^2)$ in case
of graphs, or $\Omega(n^3)$ in case of 3-CNF) the small additive
error translates into a small multiplicative error.

All the above approximation results apply to the family of so
called Constraint-Satisfaction-Problems. In some sense, these
problems can express graph properties for which one imposes
restrictions on {\bf pairs} of vertices, such as $k$-colorability.
These techniques thus fall short from applying to properties as
simple as Triangle-freeness, where the restriction is on triples
of vertices. The techniques we develop in order to obtain Theorem
\ref{t1} enable us to handle restrictions that apply to {\em
arbitrarily} large sets of vertices.

We briefly mention that $E_{{\cal P}}$ is related to packing
problems of graphs. In \cite{HR} and \cite{Yu} it was shown that
by using linear programming one can approximate the packing number
of a graph. In Section \ref{open} we explain why this technique
does not allow one to approximate $E_{{\cal P}}$.

\subsubsection{Algorithmic applications of Szemer\'edi's
Regularity Lemma} The authors of \cite{ADLRY} gave a polynomial time
algorithmic version of Szemer\'edi's Regularity Lemma. They used it
to prove that Theorem \ref{t1} holds for the $k$-colorability
property. The running time of their algorithm was improved by
Kohayakawa, R\"odl and Thoma \cite{KRT}. Frieze and Kannan \cite{FK}
further used the algorithmic version of the regularity lemma, to
obtain approximation schemes for additional graph problems.

Theorem \ref{t1} is obtained via the
algorithmic version of a strengthening of the standard regularity
lemma, which was proved in \cite{AFKS}, and it seems that these
results cannot be obtained using the standard regularity lemma.

\subsubsection{Tolerant Property-Testing} In standard
Property-Testing (see \cite{F} and \cite{Ron}) one wants to
distinguish between the graphs $G$ that satisfy a certain graph
property ${\cal P}$, or equivalently those $G$ for which $E_{{\cal
P}}(G)=0$, from those that satisfy $E_{{\cal P}}(G)
> \epsilon$. The main goal in designing property-testers is to
reduce their query-complexity, namely, minimize the number of
queries of the form "are $i$ and $j$ connected in the input
graphs?".

Parnas, Ron and Rubinfeld \cite{PRR} introduced the notion of
Tolerant Property-Testing, where one wants to distinguish between
the graphs $G$ that satisfy $E_{{\cal P}}(G) < \delta$ from those
that satisfy $E_{{\cal P}}(G) > \epsilon$, where $0 \leq \delta <
\epsilon \leq 1$ are some constants. Recently, there have been
several results in this line of work. Specifically, Fischer and
Newman \cite{FN} have recently shown that if a graph property is
testable with number of queries depending on $\epsilon$ only, then
it is also tolerantly testable for any $0 \leq \delta < \epsilon
\leq 1$ and with query complexity depending on $|\epsilon-\delta|$.
Combining this with the main result of \cite{ASmono} implies that
any monotone property is tolerantly testable for any $0 \leq \delta
< \epsilon \leq 1$ and with query complexity depending on
$|\epsilon-\delta|$. Note, that Corollary \ref{t3} implicitly states
the same. In fact, the algorithm implied by Corollary \ref{t3} is
the "natural" one, where one picks a random subset of vertices $S$,
and approximates $E_{{\cal P}}(G)$ by computing $E_{{\cal P}}$ on
the graph induced by $S$. The algorithm of \cite{FN} is far more
complicated. Furthermore, due to the nature of our algorithm if the
input graph satisfies a monotone property ${\cal P}$, namely if
$E_{{\cal P}}(G)=0$, we will always detect that this is the case.
The algorithm of \cite{FN} may declare that $E_{{\cal P}}(G)>0$ even
if $E_{{\cal P}}(G)=0$.

\subsection{Organization}

The proofs of the main results of this paper, Theorems \ref{t1}
and \ref{t4}, are independent of each other. Sections
\ref{regularity}, \ref{overview}, \ref{SecEB} and \ref{SecAlgo}
contain the proofs relevant to Theorem \ref{t1} and Sections
\ref{thm2}, \ref{St81} and \ref{St4} contain the proofs relevant
to Theorem \ref{t4}.

In Section \ref{regularity} we introduce the basic notions of
regularity and state the regularity lemmas that we use for proving
Theorem \ref{t1} and some of their standard consequences. In Section
\ref{overview} we give a high level description of the main ideas
behind our algorithms. We also state the main structural graph
theoretic lemmas, Lemmas \ref{Easy} and \ref{Big} that lie at the
core of these algorithms. The proofs of these lemmas appear in
section \ref{SecEB}. In Section \ref{SecAlgo} we give the proof of
Theorems \ref{t1} and \ref{t2} as well as a discussion about some
subtleties regarding the implementation of these algorithms.

Section \ref{thm2} contains a high-level description of the proof of
Theorem \ref{t4} as well as a description of the main tools that
we apply in this proof. In Section \ref{St81} we prove a new
Extremal Graph-Theoretic result that lies at the core of the proof
of Theorem \ref{t4}. In Section \ref{St4} we give the detailed proof
of Theorem \ref{t4}.

The final Section \ref{open} contains some concluding remarks and
open problems. Throughout the paper, whenever we relate, for
example, to a function $f_{3.1}$, we mean the function $f$ defined
in Lemma/Claim/Theorem 3.1.

\section{Regularity Lemmas and their Algorithmic Versions}\label{regularity}

In this section we discuss the basic notions of regularity, some
of the basic applications of regular partitions and state the
regularity lemmas that we use in the proof of Theorems \ref{t1}
and \ref{t2}. See \cite{KS} for a comprehensive survey on the
regularity-lemma. We start with some basic definitions. For every
two nonempty disjoint vertex sets $A$ and $B$ of a graph $G$, we
define $e(A,B)$ to be the number of edges of $G$ between $A$ and
$B$. The {\em edge density} of the pair is defined by
$d(A,B)=e(A,B)/|A||B|$.

\begin{definition}\label{regularpair}\noindent{\bf
($\gamma$-regular pair)} A pair $(A,B)$ is {\em $\gamma$-regular},
if for any two subsets $A' \subseteq A$ and $B' \subseteq B$,
satisfying $|A'| \geq \gamma|A|$ and $|B'| \geq \gamma|B|$, the
inequality $|d(A',B')-d(A,B)| \leq \gamma$ holds.
\end{definition}

Throughout the paper we will make an extensive use of the notion
of graph homomorphism which we turn to formally define.

\begin{definition}\label{homomor}\noindent{\bf
(Homomorphism)} A homomorphism from a graph $F$ to a graph $K$, is
a mapping $\varphi:V(F) \mapsto V(K)$ that maps edges to edges,
namely $(v,u) \in E(F)$ implies $(\varphi(v),\varphi(u)) \in
E(K)$.
\end{definition}

In what follows, $F \mapsto K$ denotes the fact that there is a
homomorphism from $F$ to $K$. We will also say that a graph $H$ is
homomorphic to $K$ if $H \mapsto K$. Note, that a graph $H$ is
homomorphic to a complete graph of size $k$ if and only if $H$ is
$k$-colorable.

Let $F$ be a graph on $f$ vertices and $K$ a graph on $k$ vertices,
and suppose $F \mapsto K$. Let $G$ be a graph obtained by taking a
copy of $K$, replacing every vertex with a sufficiently large
independent set, and every edge with a random bipartite graph of
edge density $d$. It is easy to show that with high probability, $G$
contains a copy of $F$ (in fact, many). The following lemma shows
that in order to infer that $G$ contains a copy of $F$, it is enough
to replace every edge with a "regular enough" pair. Intuitively, the
larger $f$ and $k$ are, and the sparser the regular pairs are, the
more regular we need each pair to be, because we need the graph to
be "closer" to a random graph. This is formulated in the lemma
below. Several versions of this lemma were previously proved in
papers using the regularity lemma (see \cite{KS}).

\begin{lemma}\label{cbmsl}
For every real $0<\eta<1$, and integers $k,f \geq 1$ there exist
$\gamma=\gamma_{\ref{cbmsl}}(\eta,k,f)$, and
$N=N_{\ref{cbmsl}}(\eta,k,f)$ with the following property. Let $F$
be any graph on $f$ vertices, and let $U_1,\ldots,U_k$ be $k$
pairwise disjoint sets of vertices in a graph $G$, where
$|U_1|=\ldots=|U_k| \geq N$. Suppose there is a mapping
$\varphi:V(F) \mapsto \{1,\ldots,k\}$ such that the following holds:
If $(i,j)$ is an edge of $F$ then $(U_{\varphi(i)},U_{\varphi(j)})$
is $\gamma$-regular with density at least $\eta$. Then
$U_1,\ldots,U_k$ span a copy of $F$.
\end{lemma}

\begin{comment}\label{GDnonIn}
Observe that the function $\gamma_{\ref{cbmsl}}(\eta,k,f)$ may and
will be assumed to be monotone non-increasing in $k$ and $f$ and
monotone non-decreasing in $\eta$. Therefore, it will be
convenient to assume that $\gamma_{\ref{cbmsl}}(\eta,k,f) \leq
\eta^2$. Similarly, we will assume that
$N_{\ref{cbmsl}}(\eta,k,f)$ is monotone non-decreasing in $k$ and
$f$. Also, for ease of future definitions (in particular those
given in (\ref{DefT})) set
$\gamma_{\ref{cbmsl}}(\eta,k,0)=N_{\ref{cbmsl}}(\eta,k,0)=1$ for
any $k \geq 1$ and $0< \eta < 1$.
\end{comment}

A partition ${\cal A}=\{V_i~|~1\leq i\leq k\}$ of the vertex set
of a graph is called an {\em equipartition} if $|V_i|$ and
$|V_{j}|$ differ by no more than $1$ for all $1\leq i < j \leq k$
(so in particular each $V_i$ has one of two possible sizes). The
{\em order} of an equipartition denotes the number of partition
classes ($k$ above). A {\em refinement} of an equipartition ${\cal
A}$ is an equipartition of the form ${\cal B}=\{V_{i,j}~|~1 \leq i
\leq k, ~1 \leq j \leq l\}$ such that $V_{i,j}$ is a subset of
$V_i$ for every $1 \leq i \leq k$ and $1 \leq j \leq l$.

\begin{definition}\label{RegPart}\noindent{\bf ($\gamma$-regular
equipartition)} An equipartition ${\cal B}=\{V_{i}~|~1 \leq i \leq
k\}$ of the vertex set of a graph is called $\gamma$-regular if
all but at most $\gamma \binom{k}{2}$ of the pairs $(V_i,V_{i'})$
are $\gamma$-regular.
\end{definition}

The Regularity Lemma of Szemer\'edi can be formulated as follows.

\begin{lemma}[\cite{Sz}]\label{SzReg}
For every $m$ and $\gamma>0$ there exists
$T=T_{\ref{SzReg}}(m,\gamma)$ with the following property: If $G$
is a graph with $n \geq T$ vertices, and ${\cal A}$ is an
equipartition of the vertex set of $G$ of order at most $m$, then
there exists a refinement ${\cal B}$ of ${\cal A}$ of order $k$,
where $m \leq k \leq T$ and ${\cal B}$ is $\gamma$-regular.
\end{lemma}

$T_{\ref{SzReg}}(m,\gamma)$ may and is assumed to be monotone
non-decreasing in $m$ and monotone non-increasing in $\gamma$.
Szemer\'edi's original proof of Lemma \ref{SzReg} was only
existential as it supplied no efficient algorithm for obtaining the
required equipartition. Alon et. al. \cite{ADLRY} were the first to
obtain a polynomial time algorithm for finding the equipartition,
whose existence is guaranteed by lemma \ref{SzReg}. The running time
of this algorithm was improved by Kohayakawa et. al. \cite{KRT} who
obtained the following result.

\begin{lemma}[\cite{KRT}]\label{SzAlg}
For every fixed $m$ and $\gamma$ there is an $O(n^2)$ time
algorithm that given an equipartition ${\cal A}$ finds
equipartition ${\cal B}$ as in Lemma \ref{SzReg}.
\end{lemma}

Our main tool in the proof of Theorem \ref{t1} is Lemma
\ref{NewReg1} below, proved in \cite{AFKS}. This lemma can be
considered a strengthening of Lemma \ref{SzReg}, as it guarantees
the existence of an equipartition and a refinement of this
equipartition that poses stronger properties compared to those
of the standard $\gamma$-regular equipartition. This stronger
notion is defined below.

\begin{definition}\label{ERegPart}\noindent{\bf (${\cal E}$-regular
equipartition)} For a function ${\cal E}(r): \mathbb{N} \mapsto
(0,1)$, a pair of equipartitions ${\cal A}=\{V_{i}~|~1 \leq i \leq
k\}$ and its refinement ${\cal B}=\{V_{i,j}~|~1 \leq i \leq k, ~1
\leq j \leq l \}$, where $V_{i,j} \subset V_i$ for all $i,j$,
are said to be ${\cal E}$-regular if
\begin{enumerate}
\item For all $1 \leq i < i' \leq k$, for all $1 \leq j,j' \leq l$
but at most ${\cal E}(k)l^2$ of them, the pair
$(V_{i,j},V_{i',j'})$ is ${\cal E}(k)$-regular.

\item All $1 \leq i < i' \leq k$ but at most ${\cal
E}(0)\binom{k}{2}$ of them are such that for all $1 \leq j,j' \leq
l$ but at most ${\cal E}(0)l^2$ of them
$|d(V_i,V_{i'})-d(V_{i,j},V_{i',j'})| < {\cal E}(0)$ holds.
\end{enumerate}
\end{definition}

It will be very important for what follows to observe that in
Definition \ref{ERegPart} we may use an arbitrary {\em function}
rather than a fixed $\gamma$ as in Definition \ref{RegPart} (such
functions will be denoted by ${\cal E}$ throughout the paper). The
following is one of the main results of \cite{AFKS}.

\begin{lemma}\noindent{\bf(\cite{AFKS})}\label{NewReg1}
For any integer $m$ and function ${\cal E}(r): \mathbb{N} \mapsto
(0,1)$ there is $S=S_{\ref{NewReg1}}(m,{\cal E})$ such that any
graph on at least $S$ vertices has an ${\cal E}$-regular
equipartition ${\cal A}$, ${\cal B}$ where $|{\cal A}|=k \geq m$ and
$|{\cal B}|=kl \leq S$.
\end{lemma}

In order to make the presentation self contained we briefly review
the proof of Lemma \ref{NewReg1}. Fix any $m$ and function ${\cal
E}$ and put $\zeta = {\cal E}(0)$. Partition $G$ into $m$ arbitrary
subsets of equal size and denote this equipartition by ${\cal A}_0$.
Put $M=m$. Iterate the following task: Apply Lemma \ref{SzReg} on
${\cal A}_{i-1}$ with $m=|{\cal A}_{i-1}|$ and $\gamma={\cal
E}(M)/M^2$ and let ${\cal A}_i$ be the refinement of ${\cal
A}_{i-1}$ returned by Lemma \ref{SzReg}. If ${\cal A}_{i-1}$ and
${\cal A}_i$ form an ${\cal E}$-regular equipartition stop,
otherwise set $M=|{\cal A}_{i-1}|$ and reiterate. It is shown is
\cite{AFKS} that after at most $100/\zeta^4$ iterations, for some $1
\leq i \leq 100/\zeta^4$ the partitions ${\cal A}_{i-1}$ and ${\cal
A}_i$ form an ${\cal E}$-regular equipartition. Moreover, detecting
an $i$ for which this holds is very easy, that is, can be done in
time $O(n^2)$ (see the proof in \cite{AFKS}). Note, that one can
thus set the integer $S_{\ref{NewReg1}}(m,{\cal E})$ to be the order
of ${\cal A}_{i}$. In particular, the following is an immediate
implication of the above discussion.

\begin{prop}\label{FuncOfEps2}
If $m$ is bounded by a function of $\epsilon$ only, then for any
${\cal E}$ the integer $S=S_{\ref{NewReg1}}(m,{\cal E})$ can be
upper bounded by a function of $\epsilon$ only.
\end{prop}

The $\epsilon$ in the above proposition will be the $\epsilon$ from
the task of approximating $E_{{\cal P}}$ within an error of
$\epsilon$ in Theorem \ref{t1}. Also, in our application of Lemma
\ref{NewReg1} the function ${\cal E}$ will (implicitly) depend on
$\epsilon$. For example, it will be convenient to set ${\cal
E}(0)=\epsilon$. However, it follows from the definition of
$S_{\ref{NewReg1}}(m,{\cal E})$ given above that even in this case
it is possible to upper bound $S_{\ref{NewReg1}}(m,{\cal E})$ by a
function of $\epsilon$ only.

In order to design our algorithm we will need to obtain the
equipartitions ${\cal A}$ and ${\cal B}$ that appear in the
statement of Lemma \ref{NewReg1}. However, note that by the overview
of the proof of Lemma \ref{NewReg1} given above, in order to obtain
this partition one can use Lemma \ref{SzAlg} as an efficient
algorithm for obtaining the regular partitions. Moreover, by
Proposition \ref{FuncOfEps2} whenever we apply either ${\cal E}$ or
Lemma \ref{SzAlg} we are guaranteed that $m$ (which in the above
overview was $M$) is upper bounded by some function of $\epsilon$
and $\gamma$ is lower bounded by some function of $\epsilon$. This
means that each of the at most $100/\zeta^4$ applications of Lemma
\ref{FuncOfEps2} takes $O(n^2)$ time. We thus get the following:

\begin{prop}\label{FuncOfEps3}
If $m$ is bounded by a function of $\epsilon$ only, then for any
${\cal E}$ there is an $O(n^2)$ algorithm for obtaining the
equipartitions ${\cal A}$ and ${\cal B}$ of Lemma \ref{NewReg1}.
\end{prop}

\section{Overview of the Proof of Theorem \ref{t1}}\label{overview}

We start with a convenient way of handling a monotone graph
property.

\begin{definition}[Forbidden Subgraphs]\label{forbidden}
For a monotone graph property $\mathcal{P}$, define
$\mathcal{F}=\mathcal{F_P}$ to be the set of graphs which are
minimal with respect to not satisfying property $\mathcal{P}$. In
other words, a graph $F$ belongs to $\mathcal{F}$ if it does not
satisfy $\mathcal{P}$, but any graph obtained from $F$ by removing
an edge or a vertex, satisfies $\mathcal{P}$.
\end{definition}
As an example of a family of forbidden subgraphs, consider $\mathcal{P}$
which is
the property of being 2-colorable. Then $\mathcal{F_P}$ is the set
of all odd-cycles. Clearly, a graph satisfies $\mathcal{P}$ if and
only it contains no member of $\mathcal{F_P}$ as a (not
necessarily induced) subgraph. We say that a graph is
$\mathcal{F}$-free if it contains no (not necessarily induced)
subgraph $F \in \mathcal{F}$. Clearly, for any family
$\mathcal{F}$, being $\mathcal{F}$-free is a monotone property.
Thus, the monotone properties are precisely the graph properties that
are equivalent to being $\mathcal{F}$-free for some family
$\mathcal{F}$. In order to simplify the notation, it will be
simpler to talk about properties of type $\mathcal{F}$-free rather
than monotone properties. To avoid confusion we will henceforth
denote by $E_{{\cal F}}(G)$ the value of $E_{{\cal P}}(G)$, where
$\mathcal{F}=\mathcal{F_P}$ as above.

The main idea we apply in order to obtain the algorithmic results of
this paper is quite simple; given a graph $G$, a family of forbidden
subgraphs ${\cal F}$ and $\epsilon > 0$ we use Lemma \ref{NewReg1}
with appropriately defined parameters in order to construct in
$O(n^2)$ time a weighted complete graph $W$, of size depending on
$\epsilon$ but {\bf independent} of the size of $G$, such that a
solution of a certain "related" problem on $W$ gives a good
approximation of $E_{{\cal F}}(G)$. As $W$ will be of size
independent of the size of $G$, we may and will use exhaustive
search in order to solve the "related" problem on $W$. In what
follows we give further details on how to define $W$ and the
"related" problem that we solve on $W$.

We start with the simplest case, where the property is that of
being triangle-free, namely ${\cal F}=\{K_3\}$. Let $W$ be some
weighted complete graph on $k$ vertices and let $0 \leq w(i,j)
\leq 1$ denote the weight of the edge connecting $i$ and $j$ in
$W$. Let $E_{{\cal F}}(W)$ be the natural extension of the
definition of $E_{{\cal F}}(G)$ to weighted graphs, namely,
instead of just counting how many edges should be removed in order
to turn $G$ into an ${\cal F}$-free graph, we ask for the edge set
of minimum weight with the above property. Let $G$ be a
$k$-partite graph on $n$ vertices with partition classes
$V_1,\ldots,V_k$ of equal size $n/k$. Suppose for every $i < j$ we
have $d(V_i,V_j)=w(i,j)$ (recall that $d(V_i,V_j)$ denotes the
edge density between $V_i$ and $V_j$). In some sense, $W$ can be
considered a weighted approximation of $G$, but to our
investigation a more important question is whether $W$ can be used
in order to estimate $E_{{\cal F}}(G)$? In other words, is it true
that $E_{{\cal F}}(G) \approx E_{{\cal F}}(W)$?

It is easy to see that $E_{{\cal F}}(G) \leq E_{{\cal F}}(W)$.
Indeed, given a set of edges $S$, whose removal turns $W$ into a
triangle free graph, we simply remove all edges connecting $V_i$ and
$V_j$ for every $(i,j) \in S$. The main question is whether the
other direction is also true. Namely, is it true that if it is
possible to remove $\alpha n^2$ from $G$ and thus make it triangle
free, then it is possible to remove from $W$ a set of edges of total
weight approximately $\alpha k^2$ and thus make it triangle-free? If
true this will mean that by computing $E_{{\cal F}}(W)$ we also
approximately compute $E_{{\cal F}}(G)$. Unfortunately, this
assertion is false in general, as the minimal number of edge
modifications that are enough to make $G$ triangle-free, may involve
removing {\em some} and not {\em all} the edges connecting a pair
$(V_i,V_j)$, and in $W$ we can remove only edges and not parts of
them. It thus seems natural to ask what kind of restrictions should
we impose on $G$ (or more precisely on the pairs $(V_i,V_j)$) such
that the above situation will be impossible, namely, that the
optimal way to turn $G$ into a triangle free graph will involve
removing either none or all the edges connecting a pair $(V_i,V_j)$
(up to some small error). This will clearly imply that we also have
$E_{{\cal F}}(G) \approx E_{{\cal F}}(W)$.

One natural restriction is that the pairs $(V_i,V_j)$ would be
random bipartite graphs. While this restriction indeed works it is
of no use for our investigation as we are trying to design an
algorithm that can handle arbitrary graphs and not necessarily
random graphs. One is thus tempted to replace random bipartite graph
with $\gamma$-regular pairs for some small enough $\gamma$.
Unfortunately, we did not manage to prove that there is a small
enough $\gamma>0$ ensuring that even if all pairs $(V_i,V_j)$ are
$\gamma$-regular then $E_{{\cal F}}(G) \approx E_{{\cal F}}(W)$. In
order to circumvent this difficulty we use the stronger notion of
${\cal E}$-regularity defined in Section \ref{regularity}. As it
turns out, if one uses an appropriately defined function ${\cal E}$,
then if all pairs $(V_i,V_j)$ are ${\cal E}(k)$-regular, one can
infer that $E_{{\cal F}}(G) \approx E_{{\cal F}}(W)$. This result is
(essentially) formulated in Lemma \ref{Easy}.

In the above discussion we considered the case ${\cal F}=\{K_3\}$.
So suppose now that ${\cal F}$ is an arbitrary (possibly infinite)
family of graph. Suppose we use a weighted complete graph $W$ on
$k$ vertices as above in order to approximate some $k$-partite
graph. The question that naturally arises at this stage is what
problem should we try to solve on $W$ in order to get an
approximation of $E_{{\cal F}}(G)$. It is easy to see that $G$ may
be very far from being ${\cal F}$-free, while at the same time $W$
can be ${\cal F}$-free, simply because ${\cal F}$ does not contain
graphs of size at most $k$. As an example, consider the case,
where the property is that of containing no copy of the complete
bipartite graph with two vertices in each side, denoted $K_{2,2}$.
Now, if $G$ is the complete bipartite graph $K_{n/2,n/2}$ then it is
very far from being $K_{2,2}$-free. However, in this case $W$ is
just an edge that spans no copy of $K_{2,2}$.

It thus seems that we must solve a {\em different} problem on $W$.
To formulate this problem we need the following definitions.

\begin{definition}\noindent{\bf(${\cal
F}$-homomorphism-free)}\label{DefHomFree} For a family of graphs
${\cal F}$, a graph $W$ is called ${\cal F}$-homomorphism-free if
$F \not \mapsto W$ for any $F \in {\cal F}$.
\end{definition}

We now define a measure analogous to $E_{{\cal F}}$ but with
respect to making a graph ${\cal F}$-homomorphism-free. Note that
we focus on weighted graphs.

\begin{definition}\noindent{\bf(${\cal H}_{\cal F}(W)$)}
For a family of graphs $\mathcal{F}$ and a weighted complete graph
$W$ on $k$ vertices, let ${\cal H'}_{\cal F}(W)$ denote the
minimum total weight of a set of edges, whose removal from $W$
turns it into an ${\cal F}$-homomorphism-free graph. Define,
${\cal H}_{\cal F}(W) = {\cal H'}_{\cal F}(W)/k^2$.
\end{definition}

Note, that in Definition \ref{DefHomFree} the graph $W$ is an
unweighed not necessarily complete graph. Also, observe that when
${\cal F}=\{K_3\}$ then we have ${\cal H}_{\cal F}(W)=E_{\cal
F}(W)$. As it turns out, the "right" problem to solve on $W$ is to
compute ${\cal H}_{\cal F}(W)$. This is formulated in the following
key lemma, whose proof appears in Section \ref{SecEB}:

\begin{lemma}\label{Easy}\noindent{\bf (The Key Lemma)}
For every family of graphs $\mathcal{F}$, there are functions
$N_{\ref{Easy}}(k,\epsilon)$ and $\gamma_{\ref{Easy}}(k,\epsilon)$
with the following property\footnote{The functions
$N_{\ref{Easy}}(k,\epsilon)$ and $\gamma_{\ref{Easy}}(k,\epsilon)$
will also (implicitly) depend on ${\cal F}$.}: Let $W$ be any
weighted complete graph on $k$ vertices and let $G$ be any
$k$-partite graph with partition classes $V_1,\ldots,V_k$ of equal
size such that
\begin{enumerate}
\item $|V_1|= \ldots =|V_k| \geq N_{\ref{Easy}}(k,\epsilon)$.

\item All pairs $(V_i,V_j)$ are
$\gamma_{\ref{Easy}}(k,\epsilon)$-regular.

\item For every $1 \leq i < j \leq k$ we have $d(V_i,V_j)=w(i,j)$.
\end{enumerate}
Then, $~ E_{\cal F}(G) \geq {\cal H}_{\cal F}(W) - \epsilon\;.$
\end{lemma}

It is easy to argue as we did above and prove that $E_{\cal F}(G)
\leq {\cal H}_{\cal F}(W)$ in Lemma \ref{Easy} (see the proof of
Lemma \ref{Big}), however we will not need this (trivial) direction.
It is important to note that while Lemma \ref{Easy} is very strong
as it allows us to approximate $E_{\cal F}(G)$ via computing ${\cal
H}_{\cal F}(W)$ (recall that $W$ is intended to be very small
compared to $G$) its main weakness is that it requires the
regularity between each of the pairs to be a function of $k$, which
denotes the number of partition classes, rather than depending
solely on the family of graphs ${\cal F}$. We note that even if
${\cal F}=\{K_3\}$ as discussed above, we can only prove Lemma
\ref{Easy} with a regularity measure that depends on $k$. This
supplies some explanation as to why Lemma \ref{SzReg} (the standard
regularity lemma) is not sufficient for our purposes; note that the
input to Lemma \ref{SzReg} is some fixed $\gamma
> 0$ and the output is a $\gamma$-regular equipartition with
number of partition classes that depends on $\gamma$ (the function
$T_{\ref{SzReg}}(m,\gamma)$). Thus, even if all pairs are
$\gamma$-regular, this $\gamma$ may be very large when considering
the number of partition classes returned by Lemma \ref{SzReg} and
the regularity measure which Lemma \ref{Easy} requires. Hence, the
standard regularity lemma cannot help us with applying Lemma
\ref{Easy}. In order to overcome this problem we use the notion of
${\cal E}$-regular partitions and the stronger regularity-lemma
given in Lemma \ref{NewReg1}, which, when appropriately used, allows
us to apply Lemma \ref{Easy} in order to obtain Lemma \ref{Big}
below, from which Theorem \ref{t1} follows quite easily. The proof
of this lemma appears in Section \ref{SecEB}.

\begin{lemma}\label{Big}
For any $\epsilon > 0$ and family of graphs ${\cal F}$ there are
functions $N_{\ref{Big}}(r)$ and ${\cal E}_{\ref{Big}}(r)$
satisfying the following\footnote{The functions $N_{\ref{Big}}(r)$
and ${\cal E}_{\ref{Big}}(r)$ will also (implicitly) depend on
$\epsilon$ and ${\cal F}$.}: Suppose a graph $G$ has an ${\cal
E}_{\ref{Big}}$-regular equipartition ${\cal A}=\{V_i~|~1 \leq i
\leq k \}$, ${\cal B}=\{V_{i,j}~|~1 \leq i \leq k,~ 1 \leq j \leq
l \}$, where
\begin{enumerate}
\item $k \geq 1/\epsilon$.

\item $|V_{i,j}| \geq N_{\ref{Big}}(k)$ for every
$1 \leq i \leq k$ and $1 \leq j \leq l$.
\end{enumerate}
Let $W$ be a weighted complete graph on $k$ vertices with
$w(i,j)=d(V_i,V_j)$. Then,
$$
|E_{\cal F}(G) - {\cal H}_{\cal F}(W)| \leq \epsilon\;.
$$
\end{lemma}

Using the algorithmic version of Lemma \ref{NewReg1}, which is given
in Proposition \ref{FuncOfEps3}, we can rephrase the above lemma in
a more algorithmic way, which is more or less the algorithm of
Theorem \ref{t1}: Given a graph $G$ we use the $O(n^2)$ time
algorithm of Proposition \ref{FuncOfEps3} in order to obtain the
equipartition described in the statement of Lemma \ref{Big}. We then
construct the graph $W$ as in Lemma \ref{Big}, and finally use
exhaustive search in order to precisely compute ${\cal H}_{\cal
F}(W)$. By Lemma \ref{Big}, this gives a good approximation of
$E_{\cal F}(G)$. The proof of Theorem \ref{t1} appears in Section
\ref{SecAlgo}.

\section{Proofs of Lemmas \ref{Easy} and \ref{Big}}\label{SecEB}

In this section we apply our new structural technique in order to
prove Lemmas \ref{Easy} and \ref{Big}. Regretfully, it is hard to
precisely state what are the ingredients of this technique. Roughly
speaking, it uses the notion of ${\cal E}$-regularity in order to
partition the edges of a graph into a bounded number of edge sets,
which have regular-partitions that are almost identical\footnote{Two
regular partitions $V_1,\ldots,V_k$ and $U_1,\ldots,U_k$ are
identical if $d(V_i,V_j)=d(U_i,U_j)$} and more importantly, the
regularity-measure of each of the bipartite graphs in each of the
edge sets can be a function of the number of clusters.

We start this section with some definitions that will be very useful
for the proof of Lemma \ref{Easy}.

\begin{definition}\label{GraphsWithSubgraph} For any (possibly
infinite) family of graphs $\mathcal{F}$, and any integer $r$ let
$\mathcal{F}_r$ be the following set of graphs: A graph $R$
belongs to $\mathcal{F}_r$ if it has at most $r$ vertices and
there is at least one $F \in \mathcal{F}$ such that $F \mapsto R$.
\end{definition}

\begin{definition}\label{DefineP} For any family of
graphs $\mathcal{F}$ and integer $r$ for which $\mathcal{F}_r \neq
\emptyset$, define
\begin{equation}\label{DefineWEq}
\Psi_{\mathcal{F}}(r)=\max_{R \in \mathcal{F}_r} ~\min_{\{F \in
\mathcal{F}: F \mapsto R\}}|V(F)|.
\end{equation}
Define $\Psi_{\mathcal{F}}(r)=0$ if $\mathcal{F}_r = \emptyset$.
Therefore, $\Psi_{\mathcal{F}}(r)$ is monotone non-decreasing in
$r$.
\end{definition}

Practicing definitions, note that if ${\cal F}$ is the family of
odd cycles, then ${\cal F}_{k}$ is precisely the family of
non-bipartite graphs of size at most $k$. Also, in this case
$\Psi_{{\cal F}}(k)=k$ when $k$ is odd, and $\Psi_{{\cal
F}}(k)=k-1$ when $k$ is even. The "right" way to think of the
function $\Psi_{{\cal F}}$ is the following: Let $R$ be a graph of
size at most $k$ and suppose we are guaranteed that there is a
graph $F' \in {\cal F}$ such that $F' \mapsto R$ (thus $R \in
{\cal F}_{k}$). Then by this information only and {\em without}
having to know the structure of $R$ itself, the definition of
$\Psi_{\mathcal{F}}$ implies that there is a graph $F \in {\cal
F}$ of size at most $\Psi_{\mathcal{F}}(k)$, such that $F \mapsto
R$.

The function $\Psi_{\mathcal{F}}$ has a critical role in the proof
of Lemma \ref{Easy}. While proving this lemma we will use Lemma
\ref{cbmsl} in order to derive that some $k$ sets of vertices, which
are regular enough, span some graph $F \in {\cal F}$. Roughly
speaking, the main difficulty will be that we will not know the size
of $F$, and as a consequence will not know the regularity measure
between these sets that is sufficient for applying Lemma \ref{cbmsl}
on these $k$ sets (this quantity is
$\gamma_{\ref{cbmsl}}(\eta,k,|V(F)|)$). However, we {\em will} know
that there is {\em some} $F' \in {\cal F}$ which is spanned by
these sets. The function $\Psi_{\mathcal{F}}(r)$ will thus be very
useful as it supplies an upper bound for the size of the smallest $F
\in {\cal F}$ which is spanned by these sets. See Proposition
\ref{SimW'}, where $\Psi_{\mathcal{F}}(r)$ has a crucial role.

\bigskip

\noindent{\bf Proof of Lemma \ref{Easy}:} Given $\epsilon$ and $k$
let
\begin{equation}\label{DefT}
T=T(k,\epsilon)=T_{\ref{SzReg}}(k,\gamma_{\ref{cbmsl}}(\epsilon/2,k,\Psi_{{\cal
F}}(k))).
\end{equation}
We prove the lemma with $\gamma_{\ref{Easy}}(k,\epsilon)$ and
$N_{\ref{Easy}}(k,\epsilon)$ satisfying

\begin{equation}\label{DefEasy1}
\gamma_{\ref{Easy}}(k,\epsilon)=\min(\epsilon/2,~ 1/T),
\end{equation}

\begin{equation}\label{DefEasy2}
N_{\ref{Easy}}(k,\epsilon)=T \cdot
N_{\ref{cbmsl}}(\epsilon/2,k,\Psi_{{\cal F}}(k))
\end{equation}

Suppose $G$ is a graph on $n$ vertices, in which case each set
$V_i$ is of size $\frac{n}{k}$. We may thus show that one must
remove at least ${\cal H}_{\cal F}(W) \cdot n^2 - \epsilon n^2$
edges from $G$ in order to make it ${\cal F}$-free. To this end,
it is enough to show that if there is a graph $G'$ that is
obtained from $G$ by removing less than ${\cal H}_{\cal F}(W)
\cdot n^2 - \epsilon n^2$ edges and spans no $F \in {\cal F}$ then
it is possible to remove from $W$ a set of edges of total weight
less than ${\cal H}_{\cal F}(W) \cdot k^2$ and obtain a graph
$W'$ that is ${\cal F}$-homomorphism-free. This will obviously
be a contradiction.

Assume such a $G'$ exists and apply Lemma \ref{SzReg} on it with
$\gamma=\gamma_{\ref{cbmsl}}(\frac12\epsilon,k,\Psi_{{\cal F}}(k))$
and $m=k$ (we use $m=k$ as $G$ is already partitioned into $k$
subsets $V_1,\ldots,V_k$). For the rest of the proof we denote by
$V_{i,1},\ldots,V_{i,l}$ the partition of $V_i$ that Lemma
\ref{SzReg} returns. Recall that as $|V_1|=\ldots=|V_k|$ and Lemma
\ref{SzReg} partitions a graph into subsets of equal size, then all
the sets $V_i$ are partitioned into the same number $l$ of subsets.
Note also that by Lemma \ref{SzReg} and the definition of $T$ in
(\ref{DefT}) we have $l < T$. Observe, that $T$ is in fact an upper
bound for the {\em total} number of partition classes $V_{i,j}$).

By Lemma \ref{SzReg} (recall that by Comment \ref{GDnonIn} we may
assume $\gamma_{\ref{cbmsl}}(\frac12\epsilon,k,\Psi_{{\cal F}}(k))
\leq \frac12\epsilon$), we are guaranteed that out of the $lk$ sets
$V_{i,j}$ at most $\frac{\epsilon}{2} \binom{lk}{2}$ pairs are not
$\gamma_{\ref{cbmsl}}(\frac12\epsilon,k,\Psi_{{\cal
F}}(k))$-regular. We define a graph $G''$, which is obtained from
$G'$ by removing all the edges connecting pairs
$(V_{i,i'},V_{j,j'})$ that are not
$\gamma_{\ref{cbmsl}}(\frac12\epsilon,k,\Psi_{{\cal
F}}(k))$-regular, and all edges connecting pairs
$(V_{i,i'},V_{j,j'})$ for which their edge density in $G'$ is
smaller than $\frac12\epsilon$.

\begin{prop}\label{SimAverage}
There are $k$ sets $V_{1,t_1},\ldots,V_{k,t_k}$ such that the
graphs induced by $G$ and $G''$ on these $k$ sets differ by less
than ${\cal H}_{\cal F}(W) \cdot \frac{n^2}{l^2}-\frac{\epsilon
n^2}{2l^2}$ edges.
\end{prop}

\noindent{\bf Proof:} We first claim that $G''$ is obtained from
$G'$ by removing less than $\frac{\epsilon}{2} n^2$ edges. To see
this note that the number of edges connecting a pair
$(V_{i,i'},V_{j,j'})$ is at most $(n/kl)^2$. As there are at most
$\frac{\epsilon}{2} \binom{lk}{2}$ pairs that are not
$\gamma_{\ref{cbmsl}}(\frac12\epsilon,k,\Psi_{{\cal F}}(k))$-regular,
we remove at most $\frac{\epsilon}{4} n^2$ edges due to such pairs.
Finally, as due to pairs, whose edge density is at most
$\frac12\epsilon$, we remove at most $\binom{kl}{2}
\frac{\epsilon}{2}(n/kl)^2 \leq \frac{\epsilon}{4} n^2$ edges, the
total number of edges removed is at most $\frac{\epsilon}{2} n^2$,
as needed.

As we assume that $G'$ is obtained from $G$ by removing less than
${\cal H}_{\cal F}(W) \cdot n^2-\epsilon n^2$ edges, we get from
the previous paragraph that $G''$ is obtained from $G$ be removing
less than ${\cal H}_{\cal F}(W) \cdot n^2-\frac{\epsilon}{2}n^2$
edges. Suppose for every $1 \leq i \leq k$ we randomly and
uniformly pick one of the sets $V_{i,1},\ldots,V_{i,l}$. The
probability that an edge, which belongs to $G$ and not to $G''$,
is spanned by these $k$ sets is $l^{-2}$. As $G$ and $G''$ differ
by less than ${\cal H}_{\cal F}(W) \cdot
n^2-\frac{\epsilon}{2}n^2$ edges, we get that the expected number
of such edges is less than ${\cal H}_{\cal F}(W) \cdot
\frac{n^2}{l^2}-\frac{\epsilon n^2}{2l^2}$ and therefore there
must be a choice of $k$ sets that span less than this number of
such edges. $\qed$

\bigskip

We are now ready to arrive at a contradiction by showing that if it
is possible to remove less than ${\cal H}_{\cal F}(W) \cdot n^2 -
\epsilon n^2$ edges from $G$ and thus turn it into an ${\cal
F}$-free graph $G'$, then we can remove from $W$ a set of edges of
total weight less than ${\cal H}_{\cal F}(W) \cdot k^2$ and thus
turn it into an ${\cal F}$-homomorphism-free graph $W'$. Let
$V_{1,i_1},\ldots,V_{k,i_k}$ be the $k$ sets satisfying the
condition of Proposition \ref{SimAverage} and obtain from $W$ a
graph $W'$ by removing from $W$ edge $(i,j)$ if and only if the
density of $(V_{i,t_i},V_{j,t_j})$ in $G''$ is 0.

\begin{prop}\label{SimW'}
$W'$ is ${\cal F}$-homomorphism-free.
\end{prop}

\noindent{\bf Proof:} Assume $F' \mapsto W'$ for some $F' \in {\cal
F}$. As $W'$ is a graph of size $k$ this means (recall Definition
\ref{DefineP}) that there is $F \in {\cal F}$ of size at most
$\Psi_{{\cal F}}(k)$ such that $F \mapsto W'$. Let $\varphi$ be a
homomorphism from $F$ to $W'$. By definition of $\varphi$, for any
$(u,v) \in E(F)$ we have $(\varphi(u),\varphi(v))$ is an edge of
$W'$. Recall that by definition of $G''$ either the density of a
pair $(V_{i,i'},V_{j,j'})$ in $G''$ is zero, or this density is at
least $\frac12\epsilon$ and the pair is
$\gamma_{\ref{cbmsl}}(\frac12\epsilon,k,\Psi_{{\cal
F}}(k))$-regular. By definition of $W'$, this means that for every
$(u,v) \in E(F)$ the pair
$(V_{\varphi(u),t_{\varphi(u)}},V_{\varphi(v),t_{\varphi(v)}})$ has
density at least $\frac{\epsilon}{2}$ in $G''$ and is
$\gamma_{\ref{cbmsl}}(\frac12\epsilon,k,\Psi_{{\cal
F}}(k))$-regular. By item 1 of the lemma we have for all $1 \leq i
\leq k$ that $|V_i| \geq N_{\ref{Easy}}(k,\epsilon)$. By our choice
in (\ref{DefEasy2}) and the fact that $l \leq T$, the sets
$V_{i,t_i}$ must therefore be of size at least
$$
|N_{\ref{Easy}}(k,\epsilon)|/l \geq |N_{\ref{Easy}}(k,\epsilon)|/T =
N_{\ref{cbmsl}}(\frac12\epsilon,k,\Psi_{{\cal F}}(k)).
$$
Hence, the sets $V_{1,t_1},\ldots,V_{k,t_k}$ satisfy all the
necessary requirements needed in order to apply Lemma \ref{cbmsl} on
them in order to deduce that they span a copy of $F$ in $G''$
(recall, that we have already argued that $|V(F)| \leq \Psi_{{\cal
F}}(k)$). This, however, is impossible, as we assumed that $G'$ was
already ${\cal F}$-free and $G''$ is a subgraph of $G'$. $\qed$

\begin{prop}\label{SimDens}
For any $i<j$ the edge densities of $(V_i,V_j)$ and
$(V_{i,t_i},V_{j,t_j})$ satisfy in $G$
$$
|d(V_i,V_j)-d(V_{i,t_i},V_{j,t_j})| \leq \frac12\epsilon.
$$
\end{prop}
\noindent{\bf Proof:} Recall that $1/l > 1/T$ and by
(\ref{DefEasy1}) we have $1/T > \gamma_{\ref{Easy}}(k,\epsilon)$. We
infer that $|V_{i,t_i}|=|V_i|/l \geq \gamma_{\ref{Easy}}(k,\epsilon)
|V_i|$. By item 2 of the lemma, each pair $(V_i,V_j)$ is
$\gamma_{\ref{Easy}}(k,\epsilon)$-regular in $G$. Hence, by
definition of a regular pair, we must have
$|d(V_i,V_j)-d(V_{i,t_i},V_{j,t_j})| \leq
\gamma_{\ref{Easy}}(k,\epsilon) \leq \frac12\epsilon$. $\qed$

\begin{prop}\label{SimWEdge}
$W'$ is obtained from $W$ by removing a set of edges of weight
less than ${\cal H}_{\cal F}(W) \cdot k^2$.
\end{prop}

\noindent{\bf Proof:} Let $S$ be the set of edges removed from $W$
and denote by $w(S)$ the total weight of edges in $S$. Let
$e(V_{i,t_i},V_{j,t_j})$ denote the number of edges connecting the
pair $(V_{i,t_i},V_{j,t_j})$ in $G$. We claim that the following
series of inequalities, which imply that $w(S) < {\cal H}_{\cal
F}(W) \cdot k^2$, hold:

\begin{eqnarray*}
{\cal H}_{\cal F}(W) \cdot
\frac{n^2}{l^2} - \frac{\epsilon n^2}{2l^2} & > & \sum_{(i,j)\in S}e(V_{i,t_i},V_{j,t_j}) \\
&\geq& \sum_{(i,j)\in
S}(w(i,j)- \frac{\epsilon}{2})\frac{n^2}{l^2k^2} \\
&\geq& \sum_{(i,j)\in S}w(i,j)\frac{n^2}{l^2k^2}-\frac{\epsilon
n^2}{2l^2} \\
&=& w(S)\frac{n^2}{l^2k^2}-\frac{\epsilon n^2}{2l^2} \;. \\
\end{eqnarray*}

Indeed, recall that by the definition of $W'$, we have $(i,j) \in S$
if and only if the density of the pair $(V_{i,i'},V_{j,j'})$ in
$G''$ is 0, which means that all the edges connecting this pair were
removed in $G''$. As by Proposition \ref{SimAverage} the total
difference between $G$ and $G''$ is less than ${\cal H}_{\cal F}(W)
\cdot \frac{n^2}{l^2} - \frac{\epsilon n^2}{2l^2}$ we infer that the
first (strict) inequality is valid. The second inequality follows
from Proposition \ref{SimDens} together with the fact that by the
condition of the lemma we have $d(V_i,V_j)=w(i,j)$. The third
inequality is due to the fact that $W$ has $k$ vertices and thus
$|S| \leq k^2$. $\qed$

\bigskip

The sought after contradiction now follows immediately from
Propositions \ref{SimW'} and \ref{SimWEdge}. This completes the
proof of the lemma. $\qed$

\bigskip

We continue with the proof of Lemma \ref{Big}.

\bigskip

\noindent{\bf Proof of Lemma \ref{Big}:}\, We prove the lemma
with:
\begin{equation}\label{EqBig1}
{\cal E}_{\ref{Big}}(r)=
\left\{%
\begin{array}{ll}
    \frac{1}{16}\epsilon^2, & \hbox{$r=0$} \\
    \min(\frac18\epsilon r^{-2},~\frac18\epsilon^2,~\gamma_{\ref{Easy}}(r,\frac18\epsilon) ), & \hbox{$r \geq 1$} \\
\end{array}%
\right.
\end{equation}
and
$$
N_{\ref{Big}}(r) = N_{\ref{Easy}}(r,\frac18\epsilon)\;.
$$

We start with showing that $E_{\cal F}(G) \leq {\cal H}_{\cal F}(W)
+ \epsilon$. Suppose $G$ is a graph of $n$ vertices, in which case
the number of edges connecting $V_i$ and $V_j$ is $w(i,j)\cdot
\frac{n^2}{k^2}$. We first remove all the edges within the sets
$V_1,\ldots,V_k$. As $k \geq 1/\epsilon$ the total number of edges
removed in this step is at most $k\binom{n/k}{2} \leq \epsilon n^2$.

Let $S$ be the set of minimal weight whose removal turns $W$ into an
${\cal F}$-homomorphism-free graph $W'$. We claim that if for every
$(i,j) \in S$ we remove all the edges connecting $V_i$ and $V_j$ the
resulting graph $G'$ spans no copy of a graph $F \in {\cal F}$.
Suppose to the contrary that $G'$ spans a copy of $F \in {\cal F}$,
and consider the mapping $\varphi:V(F) \mapsto \{1,\ldots,k\}$
that maps every vertex of $F$ that belongs to $V_j$ to $j$. As we
have removed all the edges within the sets $V_1,\ldots,V_k$ and all
edges between $V_i$ and $V_j$ for any $(i,j) \in S$ we get that
$\varphi$ is a homomorphism from $F$ to $W'$ contradicting our
choice of $S$. Finally, note that the number of edges removed in the
second step is
$$
\sum_{(i,j) \in S} w(i,j)\cdot \frac{n^2}{k^2}=n^2 \cdot {\cal
H}_{\cal F}(W)\;.
$$
Combined with the first step the total number of edges removed is
at most $n^2 \cdot {\cal H}_{\cal F}(W) + \epsilon n^2$, as
needed.

For the rest of the proof we focus on proving ${\cal H}_{\cal F}(W)
\leq E_{\cal F}(G) + \epsilon$. Let ${\cal A}$ and ${\cal B}$ be the
two equipartitions from the statement of the lemma. Suppose for
every $1 \leq i \leq k$ we randomly, uniformly and independently
pick a set $V_{i,t_i}$ out of the sets $V_{i,1},\ldots,V_{i,l}$. Let
$P$ denote the event that (i) All the pairs
$(V_{i,t_i},V_{i',t_{i'}})$ are ${\cal E}(k)$-regular. (ii) All but
at most $\frac12 \epsilon \binom{k}{2}$ of the pairs
$(V_{i,t_i},V_{i',t_{i'}})$ satisfy
$|d(V_{i,t_i},V_{i',t_{i'}})-d(V_{i},V_{i'})| \leq {\cal E}(0)$. We
need the following observations:

\begin{prop}\label{probabil}
$P$ holds with probability at least $1-\frac12\epsilon$.
\end{prop}

\noindent{\bf Proof:}\, Fix any $i<i'$. By definition of ${\cal
E}_{\ref{Big}}$ we have ${\cal E}(k) \leq \frac18\epsilon k^{-2}$,
thus by item 1 of Definition \ref{ERegPart}, the probability that
$(V_{i,t_i},V_{i',t_{i'}})$ is not ${\cal E}(k)$-regular is at
most $\frac18 \epsilon k^{-2}$. By the union bound, the
probability that one of the pairs is not ${\cal E}(k)$-regular is
at most $\binom{k}{2} \frac18 \epsilon k^{-2} \leq
\frac14\epsilon$.

Item 2 of Definition \ref{ERegPart} can be rephrased as stating that
there are at most ${\cal E}(0)\binom{k}{2}=\frac{1}{16} \epsilon^2
\binom{k}{2}$ choices of $i<i'$ for which the probability that
$|d(V_{i,t_i},V_{i',t_{i'}})-d(V_{i},V_{i'})| > {\cal
E}(0)=\frac{1}{16} \epsilon^2$ is larger than ${\cal
E}(0)=\frac{1}{16} \epsilon^2$. Thus, the expected number of $i<i'$
for which $|d(V_{i,t_i},V_{i',t_{i'}})-d(V_{i},V_{i'})|>{\cal E}(0)$
is at most $\frac{1}{16} \epsilon^2 \binom{k}{2} \cdot 1 +
\binom{k}{2} \cdot \frac{1}{16} \epsilon^2 \leq \frac18 \epsilon^2
\binom{k}{2}$. By Markov's inequality, the probability that more
than $\frac12 \epsilon \binom{k}{2}$ of $i<i'$ violate the above
inequality is at most $\frac{\epsilon}{4}$.

As properties (i) and (ii) of event $P$ each hold with probability
at least $1-\frac14\epsilon$, we get that $P$ holds with
probability at least $1-\frac12\epsilon$. $\qed$

\begin{prop}\label{distance}
Assume event $P$ holds and denote by $G'$ the subgraph of $G$
that is spanned by the sets $V_{1,t_1},\ldots,V_{k,t_k}$. Then,
$E_{\cal F}(G') \geq {\cal H}_{\cal F}(W)-\frac12\epsilon$.
\end{prop}

\noindent{\bf Proof:}\, Let $W'$ be a weighted complete graph on
$k$ vertices satisfying $w(i,i')=d(V_{i,t_i},V_{i',t_{i'}})$.
Event $P$ assumes that all the pairs $(V_{i,t_i},V_{i',t_{i'}})$
are ${\cal E}(k)$-regular. As ${\cal E}(k) \leq
\gamma_{\ref{Easy}}(k,\frac18\epsilon)$ and the lemma assumes that
$|V_{i,j}| \geq N_{\ref{Big}}(k) =
N_{\ref{Easy}}(k,\frac18\epsilon)$ we may deduce from Lemma
\ref{Easy} that
\begin{equation}\label{EqBig5}
E_{\cal F}(G') \geq {\cal H}_{\cal F}(W')-\frac{\epsilon}{8}.
\end{equation}
Now, event $P$ also assumes that all but at most $\frac{\epsilon}{2}
\binom{k}{2}$ of the pairs $i<i'$ are such that
$|d(V_i,V_{i'})-d(V_{i,t_i},V_{i',t_{i'}})| \leq {\cal E}(0) <
\frac{\epsilon}{8}$. This means that the sum of edge weights of $W'$
differs from the sum of edge weights of $W$ by at most
$\frac{\epsilon}{2}\binom{k}{2}$ due to pairs that violate the above
inequality and by at most $\binom{k}{2} \frac{\epsilon}{8}$ due to
the other pairs. This means that the sum of edge weights of $W'$
differs from that of $W$ by at most $\frac{\epsilon}{4}k^2 +
\frac{\epsilon}{16}k^2 \leq \frac{3\epsilon}{8}k^2$. This clearly
implies that
\begin{equation}\label{EqBig6}
{\cal H}_{\cal F}(W') \geq {\cal H}_{\cal F}(W) -
\frac{3\epsilon}{8}\;.
\end{equation}
The proof now follows by combining (\ref{EqBig5}) and
(\ref{EqBig6}). $\qed$

\bigskip

Let $R$ be an arbitrary set of edges whose removal from $G$ turns it
into an ${\cal F}$-free graph. Randomly and uniformly select a set
$V_{i,t_i}$ from each of the sets $V_{i,1},\ldots,V_{i,l}$, and let
$R'$ denote the set of edges of $R$ that are spanned by these $k$
sets. We claim that the following upper and lower bound on the
expected size of $R'$ hold:

\begin{eqnarray*}
\frac{1}{l^2} \cdot |R| &=& \mathbb{E}[|R'|] \\
&\geq&
\mathbb{E}[|R'|~|~P]\cdot Prob[P] \\
&\geq&
(1-\frac{\epsilon}{2}) \cdot \mathbb{E}[|R'|~|~P] \\
&\geq& (1-\frac{\epsilon}{2}) \cdot ({\cal H}_{\cal
F}(W)-\frac{\epsilon}{2}) \cdot k^2 \frac{n^2}{(kl)^2} \\
&\geq& ({\cal H}_{\cal
F}(W)-\epsilon) \cdot \frac{n^2}{l^2}\;. \\
\end{eqnarray*}

Indeed, the equality is due to the fact than an edge of $R$ has
probability precisely $1/l^2$ to be in $R'$. The second inequality
is due to Proposition \ref{probabil}, the third is due to
Proposition \ref{distance} and the last is valid because ${\cal
H}_{\cal F}(W) \leq 1$. As we thus infer that $|R| \geq {\cal
H}_{\cal F}(W) \cdot n^2 - \epsilon n^2$ for {\em arbitrary} $R$, we
get that $E_{\cal F}(G) \geq {\cal H}_{\cal F}(W) - \epsilon$, thus
completing the proof. $\qed$

\section{Proofs of Algorithmic Results}\label{SecAlgo}

The technical lemmas proved in the previous sections enabled us to
infer that certain ${\cal E}$-regular partitions may be very useful
for approximating $E_{{\cal P}}$. In this section we apply
Proposition \ref{FuncOfEps3} in order to efficiently obtain these
partitions. We first prove Theorem \ref{t1}, while overlooking some
subtle issues. We then discuss them is detail.

\bigskip

\noindent{\bf Proof of Theorem \ref{t1}:} Fix any $\epsilon > 0$
and monotone graph property ${\cal P}$. Let
$\mathcal{F}=\mathcal{F_P}$ be the family of forbidden subgraphs
of ${\cal P}$ as in Definition \ref{forbidden}. As satisfying
$\mathcal{P}$ is equivalent to being $\mathcal{F}$-free, we focus
on approximating $E_{\cal F}(G)$. Let ${\cal E}_{\ref{Big}}(r)$
and $N_{\ref{Big}}(r)$ be the appropriate function with respect to
${\cal F}$ and $\epsilon$. Put
$S(\epsilon)=S_{\ref{NewReg1}}(1/\epsilon,{\cal E}_{\ref{Big}})$
and recall that by Proposition \ref{FuncOfEps2} the integer $S$
can indeed be upper bounded by a function of $\epsilon$.

If an input graph has less than $S(\epsilon) \cdot
N_{\ref{Big}}(S(\epsilon))$ vertices we use exhaustive search in
order to precisely compute $E_{\cal F}(G)$. Assume then that $G$ has
more than $S(\epsilon) \cdot N_{\ref{Big}}(S(\epsilon))$ vertices,
and use Proposition \ref{FuncOfEps3} with $m=1/\epsilon$ and ${\cal
E}_{\ref{Big}}(r)$ as above in order to compute the equipartition
${\cal A}=\{V_i~|~1 \leq i \leq k \}$ and its refinement ${\cal
B}=\{V_{i,j}~|~1 \leq i \leq k,~ 1 \leq j \leq l \}$ satisfying the
conditions of Lemma \ref{NewReg1}. As $m$ is bounded by a function
of $\epsilon$ we get from Proposition \ref{FuncOfEps3} that this
step takes time $O(n^2)$. Also, by Lemma \ref{NewReg1} we have $kl
\leq S$, therefore, as $G$ has at least $S(\epsilon) \cdot
N_{\ref{Big}}(S(\epsilon))$ vertices each of the sets $V_{i,j}$ is
of size at least $N_{\ref{Big}}(S(\epsilon)) \geq N_{\ref{Big}}(k)$.
Let $W$ be a weighted complete graph of size $k$ where
$w(i,j)=d(V_{i},V_{j})$. Using exhaustive search, we can now
precisely compute the value of ${\cal H}_{{\cal F}}(W)$. By Lemma
\ref{Big} we may infer that $|E_{{\cal F}}(G) - {\cal H}_{{\cal
F}}(W)| \leq \epsilon$. $\qed$

\bigskip

As we have mentioned in the introduction, one should specify how the
property ${\cal P}$ is given to the algorithm. For example, ${\cal
P}$ may be an undecidable property, in which case we cannot do
anything. We thus focus on decidable graph properties. However, even
in this case we may face some unexpected problems. Note, that for a
general infinite family of graphs ${\cal F}$ it is not clear how to
compute ${\cal H}_{{\cal F}}$ in finite time. Also, returning to the
overview of the proof of Lemma \ref{NewReg1} given in Section
\ref{regularity}, note that we have implicitly assumed that one can
compute the function ${\cal E}$, as this is needed in order to
compute the parameters with which one applies Lemma
\ref{FuncOfEps2}. A close inspection of the proofs of Lemmas
\ref{Easy} and \ref{Big} reveals that computing ${\cal E}$ involves
computing the function $\Psi_{{\cal F}}$ (see (\ref{DefT}),
(\ref{DefEasy1}) and (\ref{EqBig1})). One of the main results of
\cite{ASsep} asserts that somewhat surprisingly, there is a family
of graph properties ${\cal F}$, for which the property of being
${\cal F}$-free is decidable (in fact, in $coNP$) but at the same
time $\Psi_{{\cal F}}$ is not computable. Therefore, even if we
confine ourselves to decidable graph properties we still run into
trouble.

Suppose first that $\epsilon$ is not part of the input to the
algorithm. As we have discussed in Section \ref{regularity}, in
this case all the applications of ${\cal E}_{\ref{Big}}$ are on
inputs of size depending on $\epsilon$ only, thus the algorithm
may "keep" the answers to these (finitely many) applications of
${\cal E}_{\ref{Big}}$ as part of its description. Similarly, in
this case we may need to compute ${\cal H}_{{\cal F}}$ on graphs
of size depending on $\epsilon$ only\footnote{Recall that the size
of the graph on which we compute ${\cal H}_{{\cal F}}$ is the
number of partition classes of the ${\cal E}$-regular partition,
and this number is at most $S_{\ref{NewReg1}}(m,{\cal E})$, which
is bounded by a function of $\epsilon$.}, thus the algorithm may
"keep" the answers to these (finitely many) applications of ${\cal
H}_{{\cal F}}$ as part of its description. Observe, that we don't
need to keep the answer of ${\cal H}_{{\cal F}}$ for all the
(infinite) range of edge weights. Rather, as we only need to
approximate $E_{{\cal F}}$ within an additive error of $\epsilon$,
it is enough to consider edge weights
$\{0,\epsilon,2\epsilon,3\epsilon,\ldots,1\}$.

If we want the algorithm to be able to accept $\epsilon$ as part of
the input, then we must confine ourselves to properties for which
$\Psi_{{\cal F}}$ is computable. However, as for any reasonable
graph property this function is computable, this is not a real
constraint. For example, as we have mentioned in Section
\ref{SecEB}, if ${\cal P}$ is the property of being bipartite, then
$\Psi_{{\cal F}}(k)$ is either $k$ or $k-1$. Another natural family
of properties for which $\Psi_{{\cal F}}(k)$ is computable is that
of being $H$-free for a fixed graph $H$, as in this case
$\Psi_{{\cal F}}(k) \leq |V(H)|$. By the definition of the function
${\cal E}_{\ref{Big}}$ we get that if $\Psi_{{\cal F}}$ is
computable then so is ${\cal E}_{\ref{Big}}$. It is also not
difficult to see that if $\Psi_{{\cal F}}$ is computable then so is
${\cal H}_{{\cal F}}$. Therefore, in case $\Psi_{{\cal F}}$ is
computable, there is no problem with accepting $\epsilon$ as part of
the input.

We now turn to prove Theorem \ref{t2}. We note that the above
difficulties are also relevant for Corollary \ref{t3}, which applies
Theorem \ref{t2}, but we refrain from discussing them again.

\bigskip

\noindent{\bf Proof of Theorem \ref{t2}: (sketch)} As in the
previous proof, we focus on the property of being ${\cal F}$-free,
where ${\cal F}$ is the family of forbidden subgraphs of ${\cal P}$.
Suppose, as in the previous proof, that $G$ is a large enough graph
(in terms of $\epsilon$) as otherwise we can take $D$ to be the
entire vertex set of $G$. Assume, we {\em implicitly} apply Lemma
\ref{NewReg1} on $G$ and let ${\cal A}=\{V_i~|~1\leq i\leq k\}$,
${\cal B}=\{V_{i,j}~|~1 \leq i \leq k, ~1 \leq j \leq l\}$ be the
equipartitions returned by the lemma. Let $W$ be a weighted complete
graph on $k$ vertices, where $w(i,j)=d(V_i,V_j)$. By Lemma \ref{Big}
we have

\begin{equation}\label{Sam1}
|E_{{\cal F}}(G) - {\cal H}_{{\cal F}}(W)| \leq \epsilon\;.
\end{equation}

Let $D$ be a random set of vertices and for $1 \leq i \leq k$ let
$U_i$ denote the vertices of $D$ that belong to $V_i$, and for $1
\leq i \leq k , 1 \leq j \leq l$ let $U_{i,j}$ denote the vertices
of $D$ that belong to $V_{i,j}$. Recall that $k$ and $l$ are bounded
by functions of $\epsilon$. Using standard Chernoff Bounds (see, e.g.,
\cite{ASp}), it is easy to see that if we use a large enough sample
of vertices $D$ (but only large enough in terms of $\epsilon$), then
with high probability (whp) we will have
$|d(V_i,V_{i'})-d(U_i,U_{i'})| \leq \epsilon$ for any $i < i'$ and
$|d(V_{i,j},V_{i',j'})-d(U_{i,j},U_{i',j'})| \leq \epsilon$ for any
$i < i'$ and $j \neq j'$. Therefore, if $W'$ is a weighted complete
graph on $k$ vertices, where $w(i,j)=d(U_i,U_j)$ then

\begin{equation}\label{Sam2}
|{\cal H}_{{\cal F}}(W) - {\cal H}_{{\cal F}}(W')| \leq \epsilon\;.
\end{equation}

Furthermore, using Chernoff bounds again, one can show that whp all
the pairs $(U_i,U_{i'})$ and $(U_{i,j},U_{i',j'})$ are as regular as
$(V_i,V_{i'})$ and $(V_{i,j},V_{i',j'})$ (up to $\epsilon$).
Therefore, the graph induced by $D$, denoted $G'$, will have
equipartitions ${\cal A'},{\cal B'}$ satisfying the requirements of
Lemma \ref{NewReg1}. This means that

\begin{equation}\label{Sam3}
|E_{{\cal F}}(G') - {\cal H}_{{\cal F}}(W')| \leq \epsilon\;.
\end{equation}
As (\ref{Sam1}), (\ref{Sam2}) and (\ref{Sam3}) all hold with high
probability for any $\epsilon > 0$, we can thus make sure that
with probability at least $1-\epsilon$, we will have $|E_{{\cal
F}}(G') - E_{{\cal F}}(G)| \leq \epsilon$. This completes the
proof. $\qed$

\section{Overview of the Proof of Theorem \ref{t4}}
\label{thm2}

For the proof of Theorem \ref{t4} it will be more convenient to
denote by $E'_{{\cal P}}(G)$ the number of edge removals needed to
make $G$ satisfy ${\cal P}$, in other words $E'_{{\cal P}}(G) = n^2
\cdot E_{{\cal P}}(G)$. In particular, $E'_{H}(G)$ denotes the
number of edge removals needed to turn $G$ into an $H$-free graph.
We will also denote by $E'_{r}(G)$ the number of edge removals
needed to turn $G$ into an $r$-partite graph (or equivalently
$r$-colorable graph). Note, that approximating $E'_{{\cal P}}(G)$
within $n^{2-\delta}$ is equivalent to approximating $E_{{\cal
P}}(G)$ within $n^{-\delta}$.

The main technical result we need in order to obtain Theorem
\ref{t4} is an extension of some classical results in Extremal
Graph Theory. Recall, that Tur\'an's Theorem (see \cite{W}) states
that the largest $K_{r+1}$-free graph on $n$ vertices ($K_{r+1}$ =
complete graph on $r+1$ vertices) is precisely the largest
$r$-partite graph on $n$ vertices. Another classical result is the
Erd\H{o}s-Stone-Simonovits Theorem (see \cite{W}), which states
that for any graph $H$ of chromatic number $r+1$, the largest
$H$-free graph on $n$ vertices has at most $o(n^2)$ more edges
than the largest $r$-partite graph on $n$ vertices. As any
$r$-partite graph does not contain a copy of a graph of chromatic
number $r+1$, the above results can thus be restated as saying
that when $H=K_{r+1}$ we have $E'_H(K_n) = E'_{r}(K_n)$ and that
for any $H$ of chromatic number $r+1$ we have $E'_{r}(K_n) -o(n^2)
\leq E'_H(K_n) \leq E'_{r}(K_n)$.

The main extremal graph-theoretic tool that we use in order to
obtain Theorem \ref{t4} is the following result, which greatly
extends one of the main results of \cite{BSTT}. Note, that this
result also extends Tur\'an's Theorem and the
Erd\H{o}s-Stone-Simonovits Theorem as it states that $E'_H(G)$ and
$E'_r(G)$ are very close not only when $G$ is $K_n$ but already
when $G$ has a sufficiently large minimal degree.

\begin{theo}\label{t81}
Let $H$ be a graph of chromatic number $r+1 \geq 3$.
\begin{enumerate}
\item[(i)] If there is an edge of $H$ whose removal reduces its
chromatic number, then there is constant $\mu=\mu(H)>0$ such that
if $G=(V,E)$ is a graph on $n$ vertices of minimum degree at least
$(1-\mu)n$, then $E'_H(G) = E'_{r}(G)$.
\item[(ii)] Otherwise, there are constants $\gamma=\gamma(H)>0$
and $\mu=\mu(H)>0$ such that if $G=(V,E)$ is a graph on $n$
vertices of minimum degree at least $(1-\mu)n$, then
$$ E'_{r}(G) - O(n^{2-\gamma}) \leq E'_H(G) \leq E'_{r}(G).$$
\end{enumerate}
\end{theo}

The assertion of this theorem for the special case of $H$ being a
triangle is proved in \cite{BSTT} and in a stronger form in
\cite{BKS}. We note that the $n^{2-\gamma}$ term  in the second item
of the theorem cannot be avoided. Note, that the error term we
obtain in the second part of the theorem is better than the error
term of the classical Erd\H{o}s-Stone-Simonovits Theorem. Such
improvement of the error term was previously known (see, e.g.,
\cite{Er} and \cite{Si})
but only for the case of $G$ being $K_n$ and not for $G$ of
sufficiently high minimal degree. The proof of Theorem \ref{t81}
appears in Section \ref{St81}.

Our second tool in the proof of Theorem \ref{t4} is certain
pseudo-random graphs. An {\em $(n,d,\lambda)$-graph} is a
$d$-regular graph on $n$ vertices all of whose eigenvalues, except
the first one, are at most $\lambda$ in their absolute values.
This notation was introduced by the first author in the 80s,
motivated by the fact that if $\lambda$ is much smaller than $d$,
then such graphs have strong pseudo-random properties. In
particular, (see, e.g., \cite{ASp}, Chapter 9), in this case the
number of edges between any two sets of vertices $U$ and $W$ of
$G$ is roughly its expected value, which is $|U||W|d/n$, (see
Section \ref{St4} for the precise statement). There are many known
explicit constructions of $(n,d,\lambda)$-graphs that suffice for
our purpose here. Specifically, we can use, for example, the graph
constructed by Delsarte and Goethals and by Turyn (see
\cite{KrS}). In this graph the vertex set $V(G)$ consist of all
elements of the two dimensional vector space over $GF(q)$ ($q$ is
any prime power), so $G$ has $n=q^2$ vertices. To define the edges
of $G$ we fix a set $L$ of $k$ lines through the origin. Two
vertices $x$ and $y$ of the graph $G$ are adjacent if $x-y$ is
parallel to a line in $L$. It is easy to check that this graph is
$d=k(q-1)$-regular. Moreover, because it is a strongly regular
graph, one can compute its eigenvalues precisely and show that
besides the first one they all are either $-k$ or $q-k$.
Therefore, by choosing $k=(1-\mu)\frac{q^2}{q-1}$ we obtain an
$(n,d,\lambda)$-graph with $d=(1-\mu)n$ and $\lambda \leq
\sqrt{n}$ ($\mu$ will be chosen as the constant from Theorem
\ref{t81}).

Given a graph $F$ let $F_b$ denote the $b$-blowup of $F$, that is,
the graph obtained from $F$ by replacing every vertex $v \in V(F)$
with an independent set $I_v$, of size $b$, and by replacing every
edge $(u,v) \in E(F),$ with a complete bipartite graph, whose
partition classes are the independent sets $I_u$ and $I_v$. It is
not difficult to show (see Claim \ref{BlowColor}) that for any
integer $r$, we have $E'_r(F_b)=b^2E'_r(F)$. The final piece of
notation we need is the Boolean Or, denoted by $G_1 \cup G_2$ of two
graphs $G_1$ and $G_2$ on the same set of vertices $V$. Its set of
vertices is $V$, and its set of edges contains all edges of $G_1$
and all edges of $G_2$.

Armed with these preparations, we can now outline the proof of
Theorem \ref{t4}. Its first part is an easy application of
Tur\'an's Theorem for bipartite graphs. The proof of the second
part is more interesting. Suppose all bipartite graphs satisfy
${\cal P}$, and let $r+1 ~(~\geq 3)$ be the minimum chromatic
number of a graph that does not satisfy this property. Fix a graph
$H$ of chromatic number $r+1$ that does not satisfy ${\cal P}$ and
let $\mu$ be the constant of Theorem \ref{t81}. Consider, first,
the case $r \geq 3$. In this case we show that any efficient
algorithm that approximates $E'_{{\cal P}}(G)$ up to
$n^{2-\delta}$ will enable us to decide efficiently if a given
input graph $F=(V(F),E(F))$ is $r$-colorable. Indeed, given such
an $F$ on $m$ vertices, let $b=m^c$ where $c$ is large constant,
to be chosen appropriately. Let $F_b$ be the $b$-blowup of  $F$,
and let $F'$ be the vertex disjoint union of $r$ copies of $F_b$.
Let $G'$ be the $(n,d,\lambda)$-graph with $d=(1-\mu)n$ and
$\lambda \leq \sqrt{n}$, whose number of vertices $n$, is at least
the number of vertices of $F'$, and not more than four times of
that, and identify the vertices of $F'$ with some of those of
$G'$. Let $G=G' \cup F'$ be the Boolean Or  of these two graphs.
If $F$ is $r$-colorable, then so is its blowup $F_b$, and hence in
this case $F'$ has a proper $r$-coloring in which all color
classes have the same size. This can be extended to a partition of
the vertices of $G$ to $r$ nearly equal color classes, providing
an $r$-colorable subgraph of $G$ (which satisfies ${\cal P}$ by
our choice of $r$) that contains all edges of $F'$, and some edges
of $G'$ that do not belong to $F'$. The pseudo-random properties
of $G'$ enable us to approximate this number well.

On the other hand, if $F$ is not $r$-colorable, then any
$r$-colorable subgraph of $G$ misses at least $b^2r$ edges of
$F'$, and, by the pseudo-random properties of $G'$ cannot contain
too many edges of this graph that do not belong to $F'$. With the
right choice of $c$, this will ensure that if we can approximate
the number of edges in a maximum $r$-colorable subgraph of $G$ up
to an $n^{2-\delta}$-additive error, this will enable us to know
for sure whether $F$ is $r$-colorable or not. However, by Theorem
\ref{t81}, and as the minimum degree of our graph is at least
$(1-\mu)n$, the maximum size of an $H$-free subgraph of $G$ is
very close to the maximum size of an $r$-colorable subgraph of it,
which is therefore also very close to the maximum number of edges
in a subgraph of $G$ satisfying ${\cal P}$. This implies that
approximating well this last quantity is $NP$-hard. The case $r=2$
is similar, but here we have to use that the MAX-CUT problem is
$NP$-hard. The full details appear in Section \ref{St4}.

\section{Proof of Theorem \ref{t81}}\label{St81}

Throughout this section we will assume that the number of vertices
$n$ in our graph is sufficiently large. We first prove the first
part of Theorem \ref{t81}, which is an extension of Tur\'an's
theorem. To this end, we need a result proved for $K_{r+1}$-free
graphs by Andr\'asfai, Erd\H{o}s and S\'os \cite{AES} and in a
more general form by Erd\H{o}s and Simonovits \cite{ES}.

\begin{theo}[\cite{AES},\cite{ES}]\label{Erdos-Simonovits}
Let $H$ be a fixed graph with chromatic number $r+1 \geq 3$ which
contains an edge $e$ such that $\chi(H-e)=r$. If $G$ is an $H$-free
graph of order $n$ with minimal degree
$\delta(G)>\frac{3r-4}{3r-1}n$ then $G$ is $r$-colorable.
\end{theo}

We will also need the following simple lemma.

\begin{lemma}
\label{augment-cut} Let $r \geq 2$ be an integer and suppose $G'$
is an $r$-partite subgraph of a graph $G$ (which may be empty) such
that there are $m$ edges incident to the vertices in $V(G)
\backslash V(G')$. Then $G$ has an $r$-partite subgraph of size at
least $e(G') + \frac{r-1}{r}m$.
\end{lemma}

\noindent {\bf Proof:}\ Let $(A'_1, \ldots, A'_r)$ be the
partition of $G'$. Consider an $r$-partite subgraph $\Gamma$ of
$G$ with parts $(A_1, \ldots, A_r)$ such that $A'_i \subset A_i$
for every $i$, where we place each vertex $v \in V(G) \backslash
V(G')$ in $A_i$ randomly and independently with probability $1/r$.
All edges of $G'$ are edges of $\Gamma$, and each edge incident to
a vertex in $V(G) \backslash V(G')$ appears in $\Gamma$ with
probability $\frac{r-1}{r}$. By linearity of expectation
$\mathbb{E}\big[e(\Gamma)\big] = e(G') + \frac{r-1}{r}m$, so some
$r$-partite subgraph of $G$ has at least this many edges. $\qed$

\bigskip

In particular, by taking $G'$ to be the empty graph we obtain that
every $G$ contains an $r$-partite subgraph of size at least
$\frac{r-1}{r}e(G)$.

\bigskip

\noindent{\bf Proof of Theorem \ref{t81} part (i):}\,
We prove that $E'_H(G) = E'_{r}(G)$ for all graphs $G$ on $n$ vertices
with minimum degree
$$\delta(G) \geq \left(1-\frac{3}{4(r-1)(3r-1)}\right)n+1.$$
Let $\Gamma$ be the largest (in terms of number of edges)
$r$-partite subgraph of $G$ and let $F$ be the largest $H$-free
subgraph of $G$. To prove the first part of the theorem one needs to
show that $e(F)=e(\Gamma)$. As $H$ is not $r$-colorable we trivially
have $e(F) \geq e(\Gamma)$. In the rest of the proof we establish
that $e(\Gamma) \geq e(F)$. First, note that by Lemma
\ref{augment-cut} we have
\begin{eqnarray*}
e(\Gamma) &\geq& \frac{r-1}{r}e(G)\\&=&
\frac{r-1}{r}\left(\Big(1-\frac{3}{4(r-1)(3r-1)}\Big)n+1\right)n/2\\&=&
\frac{12r^2-16r+1}{8r(3r-1)}n^2+\frac{r-1}{2r}n.
\end{eqnarray*}
If $F$ has a vertex of degree at most $\frac{3r-4}{3r-1}n$ we delete
it and continue. We construct a sequence of graphs
$F=F_n,F_{n-1},...$, where if $F_k$ has a vertex of degree $\leq
\frac{3r-4}{3r-1} k$ we delete that vertex to obtain $F_{k-1}$. Let
$F'$ be the final graph of this sequence which has $s$ vertices and
minimal degree greater than $\frac{3r-4}{3r-1}s$. Since $F'$ is
$H$-free, by Theorem \ref{Erdos-Simonovits}, it is $r$-partite.
Therefore we have that
\begin{eqnarray*}
\frac{r-1}{2r}s^2 &\geq& e(F') \geq e(F) -\frac{3r-4}{3r-1}
\left({{n+1} \choose 2}-{{s+1} \choose 2}\right)\\
&\geq& e(\Gamma)-\frac{3r-4}{2(3r-1)}(n^2-s^2)-\frac{3r-4}{2(3r-1)}n\\
&\geq& \frac{12r^2-16r+1}{8r(3r-1)}n^2
-\frac{3r-4}{2(3r-1)}(n^2-s^2).
\end{eqnarray*}
This implies that $\frac{s^2}{2r(3r-1)} \geq \frac{n^2}{8r(3r-1)}$
and so $s \geq n/2$.

Let $X$ be the set of $n-s$ vertices which we deleted, i.e.,
$X=V(G)-V(F')$. By the minimal degree assumption there are at least
$$m \geq \delta(G)|X|-{|X| \choose 2} \geq
\frac{12r^2-16r+1}{4(r-1)(3r-1)}n(n-s)+(n-s)-\frac{(n-s)^2}{2}$$
edges incident with vertices in $X$. Thus, by Lemma
\ref{augment-cut}, the size of the largest $r$-partite subgraph of
$G$ is at least
\begin{eqnarray*}
e(\Gamma)&\geq&e(F')+\frac{r-1}{r}m \geq
e(F)-\frac{3r-4}{3r-1}\left({{n+1} \choose 2}-{{s+1} \choose 2}\right)+\frac{r-1}{r}m\\
&=& e(F)-\frac{3r-4}{2(3r-1)}(n^2-s^2)-\frac{3r-4}{2(3r-1)}(n-s)+\frac{r-1}{r}m\\
&\geq& e(F)-\frac{3r-4}{2(3r-1)}(n^2-s^2)+
\frac{r-1}{r}\left(\frac{12r^2-16r+1}{4(r-1)(3r-1)}n(n-s)-\frac{(n-s)^2}{2}\right)\\
&=& e(F)+ \frac{(n-s)(2s-n)}{4r(3r-1)} \geq e(F).
\end{eqnarray*}
This implies that $e(\Gamma) \geq e(F)$ and completes the proof.
$\qed$

\bigskip

We turn to prove Theorem \ref{t81} part (ii). To this end, we first
prove the main technical result of this section, Theorem
\ref{extension} below, which is a version of Theorem
\ref{Erdos-Simonovits} that applies to arbitrary graphs $H$. We then
apply this theorem in order to prove Theorem \ref{t81} part (ii).
The reader may want to note that this application of Theorem
\ref{extension} is similar to the way we applied Theorem
\ref{Erdos-Simonovits} in order to prove the first part of Theorem
\ref{t81}.

\begin{theo}
\label{extension} Let $H$ be a fixed graph on $h$ vertices with
chromatic number $r+1\geq 3$ and let $G$ be an $H$-free graph of
order $n$ with minimum degree $\delta(G) \geq
\left(\frac{r-1}{r}-\frac{1}{3hr^2}\right)n$. Then one can delete
at most $O\big(n^{2-(r+1)/h}\big)$ edges to make $G$
$r$-colorable.
\end{theo}

\noindent {\bf Proof:}\,
First we need the following weaker bound on $E'_r(G)$.
\begin{claim}\label{claim1}
$G$ can be made $r$-partite by deleting $o(n^2)$ edges.
\end{claim}

\noindent {\bf Proof:}\, We use the Regularity Lemma given in Lemma
\ref{SzReg}. For every constant $0<\eta<\frac{1}{12hr^2}$ let
$\gamma=\gamma_{\ref{cbmsl}}(\eta,r+1,h)<\eta^2$ be sufficiently
small to guarantee that the assertion of Lemma \ref{cbmsl}
holds\footnote{Recall that by Comment \ref{GDnonIn} we may assume
that $\gamma_{\ref{cbmsl}}(\eta,r+1,h)<\eta^2$.}. Consider a
$\gamma$-regular partition $(U_1, U_2, \ldots U_k)$ of $G$. Let $G'$
be a new graph on the vertices $1 \leq i \leq k$ in which $(i,j)$ is
an edge iff $(U_i,U_j)$ is a $\gamma$-regular pair with density at
least $\eta$. Since $G$ is an $H$-free graph and $H$ is homomorphic
to $K_{r+1}$ (as $\chi(H)=r+1$), by Lemma \ref{cbmsl}, $G'$ contains
no clique of size $r+1$. Call a vertex of $G'$ {\em good} if there
are at most $\eta k$ other vertices $j$ such that the pair
$(U_i,U_j)$ is not $\gamma$-regular, otherwise call it {\em bad}.
Since the number of non-regular pairs is at most $\gamma{k \choose
2} \leq \eta^2 k^2/2$ we have that all but at most $\eta k$ vertices
are good. By definition, the degree of each good vertex in $G'$ is
at least $\left(\frac{r-1}{r}-\frac{1}{3hr^2}\right)k-2\eta k-1$,
since deletion of the edges from non-regular pairs and sparse pairs
can decrease the degree by at most $\eta k$ each and the deletion of
edges inside the sets $U_i$ can decrease it by $1$. By deleting all
bad vertices we obtain a $K_{r+1}$-free graph on at most $k$
vertices with minimal degree at least
\begin{eqnarray*}
\left(\frac{r-1}{r}-\frac{1}{3hr^2}\right)k-3\eta k-1 &\geq&
\left(\frac{r-1}{r}-\frac{2}{3hr^2}\right)k \\&\geq&
\left(\frac{r-1}{r}-\frac{1}{3r^2}\right)k \\ &>&
\frac{3r-4}{3r-1}k.
\end{eqnarray*} Therefore, by Theorem \ref{Erdos-Simonovits}, this
graph is $r$-partite. This implies that to make $G$ $r$-partite we
can delete at most $\gamma n^2 +\eta n^2+(\eta n)\cdot n +k \cdot(n
/k)^2\leq 3\eta n^2+n^2/k= o(n^2)$ edges. $\qed$

\bigskip

Consider a partition $(V_1, \ldots, V_r)$ of the vertices of $G$
into $r$ parts which maximizes the number of crossing edges
between the parts. Then for every $x \in V_i$ and $j\not =i$ the
number of neighbors of $x$ in $V_i$ is at most the number of its
neighbors in $V_j$, as otherwise by shifting $x$ to $V_j$ we
increase the number of crossing edges. By Claim \ref{claim1}, we
have that this partition satisfies that $\sum_i e(V_i)=o(n^2)$.
Call a vertex $x$ of $G$ {\em typical} if $x \in V_i$ and has at
most $n/(10hr^2)$ neighbors in $V_i$. Note that there are at most
$o(n)$ non-typical vertices in $G$ and, in particular, every part
$V_i$ contains a typical vertex. By definition, the degree of this
vertex outside $V_i$ is at least
$\left(\frac{r-1}{r}-\frac{1}{3hr^2}\right)n-\frac{n}{10hr^2}>
\left(\frac{r-1}{r}-\frac{1}{2hr^2}\right)n$ and at most
$n-|V_i|$. Therefore $|V_i| \leq (\frac{1}{r}+\frac{1}{2hr^2})n$.
Also note that the number of neighbors in $V_{i}$ of every typical
vertex $x \in V_j, j \not = i$ is at least
\begin{eqnarray}
\label{typical}
d_{V_i}(x) &\geq& d(x)-d_{V_j}(x)-(r-2)\max_k|V_k|\nonumber\\
& \geq&
\left(\frac{r-1}{r}-\frac{1}{3hr^2}\right)n-\frac{n}{10hr^2}-
(r-2)\left (\frac{1}{r}+\frac{1}{2hr^2}\right)n \nonumber\\
&>&\left(\frac{1}{r}-\frac{r-1}{2hr^2} \right)n.
\end{eqnarray}
The next claim is an immediate corollary of the above observation.

\begin{claim}
\label{claim2} Let $U$ be a subset of $V_j$ of size at least
$(\frac{1}{2r}-\frac{1}{4hr})n$ and let $y_1,\ldots, y_k$ be an
arbitrary set of $k \leq r-1$ typical vertices outside $V_j$. Then,
there are at least $\frac{n}{2r(r+1)}$ vertices in $U$, which are
adjacent to all vertices $y_i$.
\end{claim}

\noindent{\bf Proof:}\, By definition, there are at most
$|V_j|-d_{V_j}(y_i)$ non-neighbors of $y_i$ in $V_j$ and thus there
are at most that many vertices in $U$ not adjacent to $y_i$. Delete
from $U$ any vertex, which is not a neighbor of either $y_1, y_2,
\ldots, y_k$. The remaining set is adjacent to every vertex $y_i$
and has size at least
$$
|U|-\sum_{i}\big(|V_j|-d_{V_j}(y_i)\big).
$$
Since by (\ref{typical}) the degree in $V_j$ of
every typical vertex $y_i \not \in V_j$ is at least $d_{V_j}(y_i)\geq
(\frac{1}{r}-\frac{r-1}{2hr^2})n$, we obtain
that the number of common neighbors of $y_1,\ldots, y_k$ in $U$ is at least
\begin{eqnarray*}
|U|-\sum_{i}\big(|V_j|-d_{V_j}(y_i)\big)&\geq&
k\left(\frac{1}{r}-\frac{r-1}{2hr^2}\right)n -k|V_j|+|U|\\
&\geq&k\left(\frac{1}{r}-\frac{r-1}{2hr^2}\right)n -
k\left(\frac{1}{r}+\frac{1}{2hr^2}\right)n+|U|\\
&\geq& |U|-\frac{k}{2hr}n \geq \left(\frac{1}{2r}-\frac{1}{4hr}\right)n-\frac{k}{2hr}n\\
&\geq& \left(\frac{1}{2r}-\frac{k+1}{2hr}\right)n \geq
\frac{n}{2r}-\frac{n}{2h} \geq \frac{n}{2r(r+1)}.
\end{eqnarray*}
Here we used  that $k+1 \leq r$ and $h \geq r+1$. $\qed$

\begin{claim}
\label{claim3}
For every non-typical vertex  $x \in V_i$ there are at
least $\frac{n^r}{5h(3r^2)^r}$ $r$-cliques $y_1, \ldots, y_r$ such
that $y_j \in V_j$ for all $1 \leq j \leq r$ and all vertices $y_j$ are adjacent
to $x$.
\end{claim}

\noindent {\bf Proof:}\, Without loss of generality let $i=1$ and
let $x \in V_1$ be a non-typical vertex. Since for every $j\not =1$
the number of neighbors of $x$ in $V_j$ is at least as large  as the
number of its neighbors in $V_1$ we have that
\begin{eqnarray}
\label{typical2}
d_{V_j}(x) &\geq&
\frac{d_{V_j}(x)+d_{V_1}(x)}{2}\geq
\frac{1}{2}\left(\Big(\frac{r-1}{r}-\frac{1}{3hr^2}\Big)n-(r-2)
\max_i|V_i|\right)\nonumber\\
&\geq&\frac{1}{2}\left(\Big(\frac{r-1}{r}-\frac{1}{3hr^2}\Big)n-(r-2)
\Big(\frac{1}{r}+\frac{1}{2hr^2}\Big)n\right)\nonumber\\
&\geq& \left(\frac{1}{2r}-\frac{1}{4hr}\right)n.
\end{eqnarray}

To construct the $r$-cliques satisfying the assertion of the claim,
first observe, that since $x$ is non-typical it has at least
$n/(10hr^2)$ neighbors in $V_1$ and at least
$n/(10hr^2)-o(n)>n/(15hr^2)$ of these neighbors are typical. Choose
$y_1$ to be an arbitrary typical neighbor of $x$ in $V_1$ and
continue. Suppose at step $1\leq k \leq r-1$ we already have a
$k$-clique $y_1, \ldots, y_k$ such that $y_i \in V_i$ for all $i$
and all vertices $y_i$ are adjacent to $x$. Let $U_{k+1}$ be the set
of neighbors of $x$ in $V_{k+1}$. Then, by (\ref{typical2}) we have
that $|U_{k+1}|=d_{V_{k+1}}(x) \geq (\frac{1}{2r}-\frac{1}{4hr})n$
and therefore by Claim \ref{claim2} there are at least
$\frac{n}{2r(r+1)}$ common neighbors of the vertices $y_i$ in
$U_{k+1}$. Moreover, at least
$\frac{n}{2r(r+1)}-o(n)>\frac{n}{3r^2}$ of them are typical and we
can choose $y_{k+1}$ to be any of them. Therefore at the end of the
process we indeed obtained at least
$\frac{n}{15hr^2}(\frac{n}{3r^2})^{r-1}=\frac{n^r}{5h(3r^2)^r}$
$r$-cliques with the desired property. $\qed$

\begin{claim}
\label{claim4}
Each $V_i$ contains at most $O(1)$ non-typical vertices.
\end{claim}

\noindent {\bf Proof:}\, Suppose that the number of non-typical
vertices in $V_i$ is at least $5h^2(3r^2)^r$. Consider an auxiliary
bipartite graph $F$ with parts $W_1, W_2$, where $W_1$ is the set of
some $t=5h^2(3r^2)^r$ non-typical vertices in $V_i$, $W_2$ is the
family of all $n^r$ $r$-element multi-sets of $V(G)$ such that $x
\in W_1$ is adjacent to multi-set $Y$ from $W_2$ iff $Y$ is an
$r$-clique in $G$ with exactly one vertex in every $V_j$ and all
vertices of $Y$ are adjacent to $x$. By the previous claim, $F$ has
at least $e(F) \geq t\frac{n^r}{5h(3r^2)^r}=hn^r$ edges and
therefore the average degree of a vertex in $W_2$ is at least
$d_{av}=e(F)/|W_2|=e(F)/n^r \geq h$. By the convexity of the
function $f(z)={z \choose h}$, we can find $h$ vertices $x_1,
\ldots, x_h$ in $W_1$ such that the number of their common neighbors
in $W_2$ is at least
$$m \geq \frac{\sum_{Y \in W_2} {d(Y) \choose h}}{{t \choose h}}
\geq n^r \frac{{ d_{av} \choose h}}{t^h}=\Omega\big(n^r\big).$$

Thus we proved that $G$ contains $h$ vertices $X=\{x_1, \ldots,
x_h\}$ and a family of $r$-cliques $\cal C$ of size
$m=\Omega\big(n^r\big)$ such that every clique in $\cal C$ is
adjacent to all vertices in $X$. Next we need the following
well-known lemma which appears  first implicitly in Erd\H{os}
\cite{E} (see also, e.g., \cite{F1}). It states that if an
$r$-uniform hypergraph on $n$ vertices has $m=\Omega\big(n^r\big)$
edges, then it contains a complete $r$-partite $r$-uniform
hypergraph with parts of size $h$. By applying this statement to
$\cal C$, we conclude that there are $r$ disjoint set of vertices
$A_1, \ldots, A_r$ each of size $h$ such that every $r$-tuple $a_1,
\ldots, a_r$ with $a_i \in A_i$ forms a clique which is adjacent to
all vertices in $X$. The restriction of $G$ to $X, A_1, \ldots, A_r$
forms a complete $(r+1)$-partite graph with parts of size $h$ each,
which clearly contains $H$. This contradiction shows that there are
less than $5h^2(3r^2)^r=O(1)$ non-typical vertices in $V_i$ and
completes the proof of the claim. $\qed$

\bigskip

Having finished all the necessary preparations, we are now ready to
complete the proof of Theorem \ref{extension}. Let $h_1\leq h_2\leq
\ldots \leq h_{r+1}$ be the sizes of the color-classes in an $r+1$
coloring of $H$. Clearly $h_1 \leq h/(r+1)$. Without loss of
generality, suppose that $V_1$ spans at least $2hn^{2-(r+1)/h}$
edges. By the previous claim, only at most $O(n)$ of these edges are
incident to non-typical vertices. Therefore the set of typical
vertices in $V_1$ spans at least $hn^{2-(r+1)/h}$ edges. Then, by
the well known result of K\"ovari, S\'os  and Tur\'an \cite{KST}
about Tur\'an numbers of bipartite graphs, $V_1$ contains a complete
bipartite graph $H_1=K_{h_1,h_2}$ all of whose vertices are typical.
If there are at least $h_3$ typical vertices in $V_2$ which are
adjacent to all vertices of $H_1$ then we add them to $H_1$ to form
a complete $3$-partite graph $H_2$ with parts of sizes $h_1, h_2$
and $h_3$ and continue. We claim that if at step $1\leq k \leq r-1$
there is a $k+1$-partite graph $H_k \subset \cup_{i=1}^kV_i$ with
parts of sizes $h_1, \ldots, h_{k+1}$ all of whose vertices are
typical, then we can extend it to the complete $k+2$-partite graph
$H_{k+1}$ by adding $h_{k+2}$ typical vertices from $V_{k+1}$ which
are adjacent to all vertices of $H_k$. Indeed, recall that by
(\ref{typical}) the number of neighbors in $V_{k+1}$ of every
typical vertex $x \in V_i, i \not = k+1$ is at least $d_{V_{k+1}}(x)
\geq (\frac{1}{r}-\frac{r-1}{2hr^2})n$. Let $t \leq h$ be the order
of $H_{k}$. Then, as in Claim \ref{claim2} the number of vertices in
$V_{k+1}$ which are adjacent to all vertices of $H_k$ is at least
\begin{eqnarray*}
|V_{k+1}|-t\left(|V_{k+1}|-\left(\frac{1}{r}-\frac{r-1}{2hr^2}
\right)n\right) & \geq & t \left(\frac{1}{r}-\frac{r-1}{2hr^2}
\right)n-(t-1)\left
(\frac{1}{r}+\frac{1}{2hr^2}\right)n\\
&=&\frac{n}{r} - \frac{t(r-1)+t-1}{2hr^2}n\\
&\geq &\frac{n}{r}-\frac{t}{2hr}n \geq \frac{n}{r}-\frac{n}{2r}=
\frac{n}{2r}
\end{eqnarray*}
and thus at least $n/(2r)-O(1)>h_{k+2}$ of these vertices are
typical. Continuing the above process $r-1$ steps we obtain a
complete $(r+1)$-partite graph with parts of sizes $h_1, \ldots,
h_{r+1}$, which clearly contains $H$. This contradicts our
assumption that $G$ is $H$-free and shows that every $V_i$ spans
at most $O\big(n^{2-(r+1)/h}\big)$ edges. Therefore the number of
edges we need to delete to make $G$ $r$-partite is bounded by
$\sum_ie(V_i) \leq O\big(n^{2-(r+1)/h}\big)$. This completes the
proof. $\qed$

\bigskip

\noindent {\bf Proof of Theorem \ref{t81} part (ii):}\, Let $H$ be a
fixed graph on $h$ vertices with chromatic number $r+1 \geq 3$. We
show that the constants $\gamma(H)$ and $\mu(H)$ in the assertion of
the theorem can be chosen to be $(r+1)/h$ and $1/(4hr^2)$
respectively. Let $G$ be an $H$-free graph of order $n$ with minimal
degree $\delta(G)\geq (1-\frac{1}{4hr^2})n$ and let $\Gamma$ be the
largest $r$-partite subgraph of $G$ and $F$ be a largest $H$-free
subgraph of $G$. To prove the second item of the theorem it is
enough to show that $e(\Gamma) \leq e(F) \leq
e(\Gamma)+O(n^{2-(r+1)/h})$. As $H$ is not $r$-colorable we
trivially have $e(\Gamma) \leq e(F)$. In the rest of the proof we
establish that $e(F) \leq e(\Gamma)+O(n^{2-(r+1)/h})$. By Lemma
\ref{augment-cut} we have that
$$
e(\Gamma) \geq \frac{r-1}{r}e(G)=
\frac{r-1}{r}\left(1-\frac{1}{4hr^2}\right)n^2/2=
\left(\frac{r-1}{2r}-\frac{r-1}{8hr^3}\right)n^2.
$$
If $F$ has a vertex of degree at most
$(\frac{r-1}{r}-\frac{1}{3hr^2})n$ we delete it and continue. We
construct a sequence of graphs $F=F_n,F_{n-1},...$, where if $F_k$
has a vertex of degree $\leq (\frac{r-1}{r}-\frac{1}{3hr^2}) k$ we
delete that vertex to obtain $F_{k-1}$. Let $F'$ be the final graph
of this sequence which has $s$ vertices and minimal degree greater
than $(\frac{r-1}{r}-\frac{1}{3hr^2})s$ and let $\Gamma'$ be the
largest $r$-partite subgraph of $F'$. Since $F'$ is $H$-free,
Theorem \ref{extension} implies $e(F') \leq
e(\Gamma')+O(n^{2-(r+1)/h})$. Therefore we have that
\begin{eqnarray*}
\frac{r-1}{2r}s^2 +o(n^2) &\geq& e(F') \geq e(F) -
\left(\frac{r-1}{r}-\frac{1}{3hr^2}\right)
\left({{n+1} \choose 2}-{{s+1} \choose 2}\right)\\
&\geq&
e(\Gamma)-\left(\frac{r-1}{2r}-\frac{1}{6hr^2}\right)(n^2-s^2)-
O(n)\\
&\geq& \left(\frac{r-1}{2r}-\frac{r-1}{8hr^3}\right)n^2 -
\left(\frac{r-1}{2r}-\frac{1}{6hr^2}\right)(n^2-s^2)-o(n^2).
\end{eqnarray*}
This implies that
$$\frac{s^2}{6hr^2} \geq
\left(\frac{1}{6hr^2}-\frac{r-1}{8hr^3}\right)n^2 -o(n^2)>
\left(\frac{1}{6hr^2}-\frac{1}{8hr^2}\right)n^2=
\frac{n^2}{24hr^2}$$ and so $s \geq n/2$.

Let $X$ be the set of $n-s$ vertices which we deleted, i.e.,
$X=V(G)-V(F')$. By the minimal degree assumption there are at
least
\begin{eqnarray*}
m \geq \delta(G)|X|-{|X| \choose 2} &\geq&
\left(1-\frac{1}{4hr^2}\right)n(n-s)-\frac{(n-s)^2}{2}\\&=&
(n-s)\left(\Big(\frac{1}{2}-\frac{1}{4hr^2}\Big)n+\frac{s}{2}\right)
\end{eqnarray*}
edges incident with vertices in $X$. Thus, by Lemma
\ref{augment-cut}, the size of the largest $r$-partite subgraph of
$G$ is at least

\begin{eqnarray*}
e(\Gamma) &\geq&
e(\Gamma')+\frac{r-1}{r}m \geq e(F')-O\big(n^{2-(r+1)/h}\big)+ \frac{r-1}{r}m\\
\hspace{-0.1cm} &\geq& \hspace{-0.1cm}
e(F)-\left(\frac{r-1}{r}-\frac{1}{3hr^2}\right) \left({{n+1}
\choose 2}-{{s+1} \choose 2}\right)
+\frac{r-1}{r}m-O\big(n^{2-(r+1)/h}\big)\\
\hspace{-0.1cm} &\geq& \hspace{-0.1cm} e(F)-
\left(\frac{r-1}{2r}-\frac{1}{6hr^2}\right)(n^2-s^2)+
\frac{r-1}{r}m-O\big(n^{2-(r+1)/h}\big)\\
\hspace{-0.1cm} &\geq& \hspace{-0.1cm}
e(F)-\left(\frac{r-1}{2r}-\frac{1}{6hr^2}\right)(n^2-s^2)+
(n-s)\left(\Big(\frac{r-1}{2r}-\frac{r-1}{4hr^3}\Big)n+\frac{(r-1)s}{2r}\right)
-O\big(n^{2-\frac{r+1}{h}}\big)\\
\hspace{-0.1cm} &= & \hspace{-0.1cm} e(F)+
\frac{(n-s)(2s-\frac{r-3}{r}n)}{12hr^2}-O\big(n^{2-(r+1)/h}\big)
\geq e(F) -O\big(n^{2-(r+1)/h}\big). \hspace{3.2cm} \qed
\end{eqnarray*}

\section{Proof of Theorem \ref{t4}}\label{St4}

We start with the proof of the first part of Theorem \ref{t4}. If
there is a bipartite graph $H$ that does not satisfy ${\cal P}$,
then, by the known results about the Tur\'an numbers of bipartite
graphs proved in \cite{KST}, there exists a positive $\delta >0$
such that for any large $n$, any graph with $n$ vertices and at
least $n^{2-\delta}$ edges contains a copy of $H$. Thus, given a
graph $G$ on $n$ vertices, one must delete all its edges besides,
possibly, $n^{2-\delta}$ of them, to obtain a subgraph satisfying
${\cal P}$. As certainly the edgeless graph satisfies ${\cal P}$,
this provides the required approximation in this case.

The proof of the second part is more complicated, and requires all
the preparations obtained in the previous section. Suppose all
bipartite graphs satisfy ${\cal P}$, and let $r+1\geq 3$ be the
minimum chromatic number of a graph that does not satisfy this
property. Fix a graph $H$ of chromatic number $r+1$ that does not
satisfy ${\cal P}$. We will show that any efficient algorithm that
approximates $E'_{{\cal P}}(G)$ up to $n^{2-\delta}$ will enable
us to decide efficiently how many edges we need to delete from a
given input graph $F=(V(F),E(F))$ to make it $r$-partite. For
$r\geq 3$ this problem contains the $r$-colorability problem, and
for $r=2$ it is the MAX-CUT problem and therefore it is $NP$-hard
for every $r \geq 2$.

Given a graph $F$ on $m$ vertices such that we need to delete
$\ell$ edges to make it $r$-partite, let $b=m^c$ where $c$ is a
large constant, to be chosen later. Let $F_b$ be the $b$-blowup of
$F$, and let $F'$ be the vertex disjoint union of $r$ copies of
$F_b$. Let $\mu=\mu(H)$ be the constant from Theorem \ref{t81} and
let $G'$ be the $(n,d,\lambda)$-graph with $d=(1-\mu)n$ and
$\lambda \leq \sqrt{n}$, described in Section \ref{thm2}. As the
integer $q$ in the construction discussed in Section \ref{thm2}
can be a prime power, we can always choose the number of vertices
of $G'$, which is $q^2$, to be at least the number of vertices of
$F'$, and not more than 4 times of that. In particular, we have
$n=\Theta(rmb)=\Theta\big(m^{c+1}\big)$. Identify the vertices of
$F'$ with some of those of $G'$. Let $G=G' \cup F'$ be the Boolean
Or of these two graphs.

Suppose, that instead of adding to $F'$ a pseudo-random graph
$G'$, we would put any non-edge of $F'$ in $G$ with probability
$1-\mu$. It is easy to see that in this case the expected number
of edges, which would be spanned by a set of $a$ vertices that
span $t$ edges in $F'$, would be $(1-\mu)\binom{a}{2}+\mu t$. The
following claim establishes that this is approximately what we
find when we add to $F'$ a pseudo-random graph. We then use this
claim to show that we can also estimate $E'_r(G)$ as a function of
$\ell=E'_r(F)$.

\begin{claim}\label{EdgeDist}
Let $A$ be a subset of the vertices of $G$ of size $a$ which
contains precisely $t$ edges of $F'$. Then the number of edges of
$G$ in $A$ satisfies
$$(1-\mu)\frac{a^2}{2}+\mu t-O\big(m^2n^{3/2}\big) \leq
e_G(A) \leq (1-\mu)\frac{a^2}{2}+\mu t+O\big(m^2n^{3/2}\big).$$
\end{claim}

\noindent {\bf Proof:}\, By construction, the edges of the
subgraph of $F'$ induced on the set $A$ form an edge disjoint
union of complete bipartite graphs we denote by $\Gamma_i=(U_i,
W_i), 1 \leq i \leq k$. Thus $\sum_i|U_i|W_i|=t$ and the fact that
$F'$ is a blowup of $r$ disjoint copies of $F$, which altogether
have $rm$ vertices and at most $r{m \choose 2}$ edges, implies
that $k \leq r{m \choose 2}<r m^2$. The number of edges of $G$
spanned on $A$ is  the number of edges of $G'$ inside $A$, minus
the number of edges of $G'$ spanned by the pairs $(U_i, W_i)$,
plus the number of edges of $F'$ inside $A$. To estimate this
quantity, we need the well-known fact (see, e.g, Chapter 9 of
\cite{ASp}), that the number of edges between two subsets $X,Y$ of
an $(n,d,\lambda)$-graph $G'$ satisfies
$$\left|\,e(X,Y)-\frac{|X||Y|d}{n}\,\right|\le \lambda\sqrt{|X||Y|}$$
and the fact that in such a graph
$\big|e(X)-\frac{d|X|^2}{2n}\big| \le \lambda|X|$. Therefore we
obtain that

\begin{eqnarray*}
e_{G}(A) &=& e_{G'}(A)-\sum_{i=1}^k e_{G'}(U_i,W_i)+t=
e_{G'}(A)+\sum_{i=1}^k\Big(|U_i|W_i|-e_{G'}(U_i,W_i)\Big)\\
&\geq& \frac{d|A|^2}{2n}-\lambda|A|+
\sum_{i=1}^k\left(|U_i|W_i|-\frac{d}{n}|U_i|W_i|-\lambda\sqrt{|U_i|W_i|}\right)\\
&\geq& \frac{d|A|^2}{2n}-\lambda n+\sum_{i=1}^k\big(
\mu|U_i|W_i|-\lambda n\big)
\\&=&(1-\mu)\frac{a^2}{2}+\mu \sum_{i=1}^k|U_i|W_i| -(k+1)\lambda n\\
&=& (1-\mu)\frac{a^2}{2}+\mu t-O\big(m^2n^{3/2}\big).
\end{eqnarray*}
The upper bound $e_{G}(A) \leq (1-\mu)\frac{a^2}{2}+\mu
t+O\big(m^2n^{3/2}\big)$ can be obtained similarly. $\qed$

\bigskip

Recall that the $b$-blowup $F_b$ of a graph $F$, defined in Section
\ref{thm2}, is the graph obtained from $F$ by replacing every vertex
$v \in V(F)$ with an independent set $I_v$, of size $b$, and by
replacing every edge $(u,v) \in E(F),$ with a complete bipartite
graph, whose partition classes are the independent sets $I_u$ and
$I_v$.

\begin{claim}\label{BlowColor}
For any graph $F$ and any integer $b$, we have
$E'_r(F_b)=b^2E'_r(F)$.
\end{claim}

\noindent{\bf Proof:}\, We start by showing that $E'_r(F_b) \leq
b^2E'_r(F)$. Suppose $S$ is a set of $E'_r(F)$ edges whose removal
turns $F$ into an $r$-colorable graph $F'$. Suppose we remove from
$F_b$ all the edges connecting $I_u$ and $I_v$ for any $(u,v) \in
S$. Note, that we thus remove $b^2E'_r(F)$ edges from $F_b$. We
claim that the resulting graph $F'_b$ is $r$-colorable. Indeed, let
$c:V(F) \mapsto \{1,\ldots,r\}$ be a $r$-coloring of $F'$ and note
that by definition of $F'_b$, if we color all the vertices of $I_v$
with the color $c(v)$, we get a legal $r$-coloring of $F'$.
Therefore $E'_r(F_b) \leq b^2E'_r(F)$.

To see that $E'_r(F_b) \geq b^2E'_r(F)$, let $S$ be a set of edges
whose removal turns $F_b$ into an $r$-colorable graph, and suppose
for every $v \in V(F)$ we randomly pick a single vertex from each of
the sets $I_v$. For every edge of $S$, the probability that we
picked both of its endpoints is $b^{-2}$, therefore the expected
number of edges spanned by these vertices is $|S|/b^2$. As the
removal of the edges of $S$ makes $F_b$ $r$-colorable, this in
particular applies to all of its subgraphs. Note, that for any
choice of a single vertex from each of the independent sets $I_v$,
the graph they span is isomorphic to $F$. Thus, any such choice
spans at least $E'_r(F)$ of the edges of $S$. It thus must be the
case that $|S|/b^2 \geq E'_r(F_b)$, and the proof is complete.
$\qed$

\begin{claim}\label{EstE}
The graph $G$ satisfies
\begin{equation}\label{SimRem}
\left|E'_r(G)-\Big((1-\mu)\frac{n^2}{2r}+\mu r\ell b^2\Big)\right| \leq
O(m^2n^3).
\end{equation}
\end{claim}

\noindent {\bf Proof:}\, Fix a partition of $F$ into $r$ parts which
misses exactly $\ell$ edges and consider $r$ disjoint copies of $F$.
By taking appropriately different parts in every copy of $F$ we can
partition this new graph into $r$ equal parts such that exactly
$r\ell$ edges are non-crossing. Since $F'$ is a $b$-blowup of $r$
disjoint copies of $F$, this gives a partition of $F'$ into equal
parts which misses $r\ell b^2$ edges. We can extend this to a
partition of $G$ into $r$ nearly equal sets $V(G)=V_1 \cup \ldots
\cup V_r$ which misses exactly $r\ell b^2$ edges of $F'$. Let $t_i$
be the number of edges of $F'$ inside $V_i$, then $\sum_i t_i =
r\ell b^2$. This, together with Claim \ref{EdgeDist}, implies that
it is enough to delete at most
\begin{eqnarray*} \sum_{i=1}^r e_{G}(V_i) &\leq&
\sum_{i=1}^r\left((1-\mu)\frac{|V_i|^2}{2}+\mu t_i+O\big(m^2n^{3/2}\big)\right)\\
&\leq& (1-\mu)r \frac{(n/r+1)^2}{2} + \mu \sum_{i=1}^rt_i+O\big(m^2n^{3/2}\big)\\
&=&(1-\mu)\frac{n^2}{2r}+\mu r\ell b^2+O\big(m^2n^{3/2}\big).
\end{eqnarray*}
edges to make $G$ $r$-partite and hence to satisfy property $\cal
P$.

On the other hand, by Claim \ref{BlowColor}, any partition of
$F'$, which is $b$-blowup of $r$ disjoint copies of $F$, into $r$
parts misses at least $r \ell b^2$ edges. Therefore for every
partition of the vertices of $G$ into $r$ sets there are at least
$r\ell b^2$ edges of $F'$ which are non-crossing. Let $V_1\cup
\ldots \cup V_r$ be a partition of $V(G)$ that maximizes the
number of crossing edges and let again $t_i$ be the number of
edges of $F'$ inside $V_i$ (note that in this case the sets $V_i$
are not necessarily of the same size). Using Claim \ref{EdgeDist},
together with the fact that $\sum_i t_i \geq r \ell b^2$ and the
Cauchy-Schwartz inequality, we conclude that

\begin{eqnarray*}
\sum_{i=1}^re_G(V_i) &\geq&
\sum_{i=1}^r\left((1-\mu)\frac{|V_i|^2}{2}+\mu t_i-O\big(m^2n^{3/2}\big)\right)\\
&\geq& \frac{1-\mu}{2}r \left(\frac{\sum_i|V_i|}{r}\right)^2 + \mu
r\ell
b^2-O\big(m^2n^{3/2}\big)\\
&=&(1-\mu)\frac{n^2}{2r}+\mu r\ell b^2-O\big(m^2n^{3/2}\big).
\end{eqnarray*}
This completes the proof of the claim. $\qed$

\bigskip

We are now ready to complete the proof of Theorem \ref{t4}. Choose
the constant $c$ to be sufficiently large so that $2/(c+1) <
\min(\delta,\gamma, 1/4)$. Recall, that as we chose $b=m^c$ and
$n=\Theta(m^{c+1})$, we have
\begin{equation}\label{SimRem2}
n^{2-\delta}=o(b^2),~~~~~ n^{2-\gamma}=o(b^2),~~~~~
m^2n^{3/2}=o(b^2).
\end{equation}
Also, as $G$ has minimum degree $(1-\mu)n$ we get from Theorem
\ref{t81}, that
\begin{equation}\label{SimRem3}
E'_H(G) \geq E'_r(G) - O(n^{2-\gamma}).
\end{equation}
As $H$ does not satisfy ${\cal P}$ we clearly have $E'_{{\cal
P}}(G) \geq E'_{H}(G)$. Combining this with (\ref{SimRem}),
(\ref{SimRem2}) and (\ref{SimRem3}) we get
\begin{eqnarray*}
E'_{{\cal P}}(G) \geq E'_{H}(G) & \geq & E'_r(G)- O(n^{2-\gamma})
\\ & \geq & (1-\mu)\frac{n^2}{2r}+\mu r \ell
b^2-O\big(m^2n^{3/2}\big)-
O\big(n^{2-\gamma}\big) \\
&\geq& (1-\mu)\frac{n^2}{2r}+\mu r \ell b^2-o\big(b^2\big) \;. \\
\end{eqnarray*}
Furthermore, by our choice of $r$, we get that any $r$-colorable
graph satisfies ${\cal P}$, hence we infer from (\ref{SimRem}) and
(\ref{SimRem2}) that

\begin{eqnarray*}
E'_{{\cal P}}(G) \leq E'_{r}(G) & \leq & (1-\mu)\frac{n^2}{2r}+\mu
r
\ell b^2+O\big(m^2n^{3/2}\big) \\
&\leq& (1-\mu)\frac{n^2}{2r}+\mu r \ell b^2+o\big(b^2\big) \;. \\
\end{eqnarray*}

We thus conclude that $|E'_{{\cal P}}(G) -
((1-\mu)\frac{n^2}{2r}+\mu r \ell b^2)| \leq o(b^2)$. Therefore,
if one can approximate $E'_{{\cal P}}(G)$ in time polynomial in
$n$ (and hence also in $m$) within an additive error of
$n^{2-\delta}=o(b^2)$ then one thus efficiently computes an
integer $L$, which is within an additive error of $o(b^2)$ from
$(1-\mu)\frac{n^2}{2r}+\mu r \ell b^2$. But as in this case $\ell$
is precisely the nearest integer to $(L-(1-\mu)\frac{n^2}{2r})/\mu
r b^2$, this implies that we can {\em precisely} compute the
number of edge removals, needed in order to turn the input graph
$F$ into an $r$-partite graph. This implies that the problem of
approximating $E'_{{\cal P}}(G)$ within $n^{2-\delta}$ is
$NP$-hard, and completes the proof of Theorem \ref{t4}.

\section{Concluding Remarks and Open Problems}\label{open}

\begin{itemize}

\item We have shown that for any monotone graph property ${\cal
P}$ and any $\epsilon >0$ one can approximate efficiently the
minimum number of edges that have to be deleted from an $n$-vertex
input graph to get a graph that satisfies ${\cal P}$, up to an
additive error of $\epsilon n^2$. Moreover, for any {\em dense}
monotone property, that is, a property for which there are graphs
on $n$ vertices with $\Omega (n^2)$ edges that satisfy it, it is
$NP$-hard to approximate this minimum up to an additive error of
$n^{2-\delta}$. It will be interesting to obtain similar sharp
results for the case of sparse monotone properties. In some of
these cases (like the property of containing no cycle, or the
property of containing no vertex of degree at least $2$) the above
minimum can be computed precisely in polynomial time, and in some
other cases, a few of which are treated in \cite{Asano},
\cite{AH}, \cite{Y}, a precise computation is known to be hard.
Obtaining sharp estimates for the best approximation achievable
efficiently seems difficult.

\item As we have mentioned in Section \ref{intro}, a special case of
Theorem \ref{t4} implies that for any non-bipartite $H$, computing
the smallest number of edge removals that are needed to make a graph
$H$-free is $NP$-hard. This is clearly not the case for some
bipartite graphs such as a single edge or any star. It will be interesting to
classify the bipartite graphs for which this problem is $NP$-hard.

\item It seems interesting to decide if one can obtain
a result analogous to Theorem \ref{t4} for the family of
hereditary properties.

\item A weaker version of Theorem \ref{t1} can be derived by
combining the results of \cite{ASmono} and \cite{FN}.
However, this only enables one to
approximate $E_{{\cal P}}(G)$ within an additive error $\epsilon$
in time $n^{f(\epsilon)}$, while the running time of our algorithm
is of type $f(\epsilon)n^2$.

\item
Recall that $E'_{{\cal F}}(G)$ denotes the smallest number of
edge deletions that are needed in order to make
$G$ ${\cal F}$-free. For a family of graphs
${\cal F}$, let $\nu_{{\cal F}}(G)$ denote the ${\cal F}$-packing
number of $G$, which is the size of the largest family of
edge-disjoint copies of members of ${\cal F}$, which is spanned by
$G$. Let $\nu^*_{{\cal F}}(G)$ denote the natural Linear Programming
relaxation of $\nu_{{\cal F}}(G)$. Haxell and R\"odl \cite{HR} and
Yuster \cite{Yu} have shown that $\nu_{{\cal F}}(G) \leq
\nu^*_{{\cal F}}(G) \leq \nu_{{\cal F}}(G) + \epsilon n^2$ for any
${\cal F}$ and any $\epsilon > 0$, implying that for any finite
${\cal F}~$, $\nu_{{\cal F}}(G)$
can be approximated within any additive error of $\epsilon n^2$ by
solving the Linear Program for computing $\nu^*_{{\cal F}}(G)$. One may wonder
whether it is possible to obtain Theorem \ref{t1} by
solving the natural Linear Programming relaxation of $E'_{{\cal
F}}(G)$, which we denote by $E^*_{{\cal F}}(G)$. Regretfully, this
is not the case. Linear Programming duality implies that $E^*_{{\cal
F}}(G) = \nu^*_{{\cal F}}(G)$ and by the results of \cite{HR} and
\cite{Yu} we thus have

\begin{equation}\label{HRY}
\nu_{{\cal F}}(G) \leq E^*_{{\cal F}}(G) \leq \nu_{{\cal F}}(G) +
\epsilon n^2 \;.
\end{equation}

Consider now any ${\cal F}$, which does not contain the single edge
graph and note that we trivially have $\nu_{{\cal F}}(K_n) \leq
\frac12\binom{n}{2} \leq \frac14 n^2$ (we denote by $K_n$ the
$n$-vertex complete graph). If ${\cal F}$ contains a bipartite graph
then by the theorem of K\"ovari, S\'os and Tur\'an (see Section
\ref{thm2}) we have
$E'_{{\cal F}}(K_n) > {n \choose 2}-n^{2-\delta} \geq (\frac12-o(1)) n^2$.
If on the other hand all the graphs in ${\cal F}$ are of chromatic
number $r \geq 3$ then clearly they all must contain at least
$\binom{r}{2}$ edges, and therefore we must have $\nu_{{\cal
F}}(K_n) \leq \binom{n}{2}/\binom{r}{2} \leq
\frac{n^2}{r(r-1)}$.
On the other hand, by the theorem of Erd\H{o}s-Stone-Simonovits
(see Section \ref{thm2}) $E'_{{\cal F}}(K_n) > \frac{n^2}{2(r-1)}
-o(n^2)$. In any case, we have that $\nu_{{\cal F}}(K_n)+ \delta n^2
\leq E'_{{\cal F}}(K_n)$ for some fixed $\delta=\delta({\cal F}) >
0$. Combined with (\ref{HRY}) we get that for any ${\cal F}$ not
containing the single edge graph $E^*_{{\cal F}}(K_n) + \delta n^2 <
E'_{{\cal F}}(K_n)$. Thus, the (trivial) case in which ${\cal F}$
contains a single edge is the only one for which computing
$E^*_{{\cal F}}(G)$ is guaranteed to approximate $E'_{{\cal F}}(G)$
within $\epsilon n^2$ for any $\epsilon > 0$. In fact, in this degenerate
case we actually have $E^*_{{\cal
F}}(G) = E'_{{\cal F}}(G)$.

\end{itemize}

\end{document}